\documentclass{amsart}

\usepackage{amsmath,amsthm,amssymb}

\renewcommand{\frak}{\mathfrak}
\renewcommand{\geq}{\geqslant}
\renewcommand{\leq}{\leqslant}
\renewcommand{\div}{\operatorname{div}}
\newcommand{\bbold}{\mathbb}
\newcommand{\rom}{\textup}
\newcommand{\Nil}{\operatorname{Nil}}

\newcommand{\id}{\operatorname{id}}

\newcommand{\abs}[1]{\lvert#1\rvert}

\def \Q {{\bbold Q}}
\def \Z {{\bbold Z}}
\def \N {{\bbold N}}

\def \F {{\bbold F}}
\def \m {{\frak m}}
 
\def \bar {\overline}
\def \<{\langle}
\def \>{\rangle}
\def \tilde {\widetilde}
\def \iso {\cong}

\def \einb{\hookrightarrow}

\def \rad{\operatorname{rad}}
\def \Spec{\operatorname{Spec}}
\def \alg{{\operatorname{alg}}} 
\def \sep{{\operatorname{sep}}} 
\def \pow{{\operatorname{pow}}}

\def \Frac {\operatorname{Frac}}

\def \Sol {\operatorname{Sol}}
\def \cal {\mathcal}
\def \tr {\operatorname{tr}}

\newtheorem{theorem}{Theorem}[section]

\newtheorem*{theoremA}{Theorem A}
\newtheorem*{theoremB}{Theorem B}
\newtheorem*{theoremC}{Theorem C}
\newtheorem*{theoremD}{Theorem D}
\newtheorem{lemma}[theorem]{Lemma}
\newtheorem{prop}[theorem]{Proposition}
\newtheorem{cor}[theorem]{Corollary}

\theoremstyle{definition}
\newtheorem{definition}[theorem]{Definition}
 
\theoremstyle{remark}

\newtheorem*{example}{Example}
\newtheorem*{examples}{Examples}
\newtheorem*{remarks}{Remarks}
\newtheorem*{remark}{Remark}
\newtheorem*{question}{Question}

\numberwithin{equation}{section}

\begin{document}

\title{Bounds and definability in polynomial rings}

\author{Matthias Aschenbrenner}
%\address{Mathematical Sciences Research Institute \\
%17 Gauss Way\\
%Berkeley, CA 94720, U.S.A. }

\address{Department of Mathematics\\
University of California at Berkeley\\
Evans Hall \\
Berkeley, CA 94720, U.S.A. }
\email{maschenb@math.berkeley.edu}

\thanks{Partially supported by the Mathematical Sciences Research
Institute, Berkeley, CA}

\date{June 2003}

\begin{abstract} 
We study questions around the existence of bounds and the
dependence on parameters for linear-algebraic problems in
polynomial rings over rings
of an arithmetic flavor.
In particular, we show that the module of syzygies of
polynomials $f_1,\dots,f_n\in R[X_1,\dots,X_N]$ with coefficients in 
a Pr\"ufer domain $R$ can be generated by elements whose degrees are 
bounded by a number only depending on $N$, $n$ and the degree of the $f_j$.
This implies that if $R$ is a B\'ezout domain, then
the generators can be parametrized in terms of the
coefficients of $f_1,\dots,f_n$ using
the ring operations and a certain division function, uniformly in
$R$.
\end{abstract}


\maketitle

\section*{Introduction}

The main theme of this article is the existence of bounds 
for basic operations of linear algebra in polynomial rings over (commutative)
rings of an arithmetic nature.  
The following result, shown in Section~\ref{Homogeneous-Section} 
below, is typical.

\begin{theoremA}
Given integers $N,d,n\geq 0$ there exists an integer $\beta=\beta(N,d,n)$ with
the following property: for every Pr\"ufer domain $R$ and polynomials
$f_1,\dots,f_n\in R[X]=R[X_1,\dots,X_N]$ of \rom{(}total\rom{)} degree
$\leq d$, the $R[X]$-submodule of $R[X]^n$ consisting of all solutions to
the linear homogeneous equation
$$f_1y_1+\cdots+f_ny_n=0$$
can be generated by \rom{(}finitely many\rom{)}
solutions whose components have degree $\leq\beta$.
\end{theoremA}

A classical theorem due to G.~Hermann \cite{Hermann} states that Theorem~A is
true if we replace ``Pr\"ufer domain'' by ``field''. In this case,
it is easy to see that $\beta$ can be chosen independent of $n$;
Seidenberg \cite{Seidenberg1} computed an explicit (doubly 
exponential) bound $\beta$. In \cite{maschenb-ideal2} we extended 
Hermann's result to the class of {\em almost Dedekind domains}\/
(i.e., domains all of whose localizations at maximal ideals are discrete
valuation rings) and obtained the bound $$\beta(N,d)=(2d)^{2^{O(N^2)}}.$$
In contrast to \cite{maschenb-ideal2}, 
the methods employed to prove Theorem~A in the present paper are 
rather non-constructive. They are inspired by the model-theoretic
approach (see, e.g., \cite{vdDries-Schmidt}) 
to establish the existence of uniform bounds, a difference being
our use of {\em direct products}\/ rather than {\em ultraproducts}\/
(or other versions of the Compactness Theorem of first-order logic).
We also work in the more general setting of
{\em semihereditary rings}\/ and rely in an essential way
on a theorem of Vasconcelos (see Theorem~\ref{Vasconcelos} below)
on the coherence of polynomial rings over semihereditary rings.
Theorem~A remains true for certain possibly non-reduced rings as well,
in particular for Artinian local rings of fixed length.
(See Corollary~\ref{Artinian-Cor}.)

The following theorem 
shows that the analogue of Theorem~A for 
inhomogeneous linear equations holds only in a very restricted setting:

\begin{theoremB}
For a ring $R$, the following statements are equivalent:
\begin{enumerate}
\item The 
nilradical $$\Nil(R)=\{r\in R:\text{$r^n=0$ for some $n\geq 1$}\}$$
of $R$ is nilpotent, and $R/\Nil(R)$ is von Neumann regular.
\item For all integers $N,d,n\geq 0$ there exists an integer
$\beta=\beta(N,d,n)$ with the following property:
if $f_1,\dots,f_n\in R[X]=R[X_1,\dots,X_N]$ are of degree at most $d$
such that 
$$1=f_1g_1+\cdots+f_ng_n$$
for some $g_1,\dots,g_n\in R[X]$, then there exist such $g_j$ of
degree at most $\beta$.
\end{enumerate}
\end{theoremB}

We prove Theorem~B in Section~\ref{Inhomogeneous-Section}, using
\cite{Kollar} for the implication (1)~$\Rightarrow$~(2) and
by combining
a result of Sabbagh \cite{Sabbagh} with an elementary lemma of Cherlin
\cite{Cherlin} for the converse. 
Note that condition (1) in Theorem~B is satisfied if $R$ is an
Artinian local ring, yielding a result on 
uniform bounds for the ideal membership problem
over $R$ originally proved by Schoutens \cite{Schoutens}.
Condition (1) is clearly also satisfied if $R$ is a field. This case of
our theorem is again due to Hermann \cite{Hermann}. Here
$\beta$ does not depend on the particular field. The existence
of such a bound is equivalent to the following: if
$f_0,f_1,\dots,f_n\in\Z[C,X]$ ($C=(C_1,\dots,C_M)$ being parametric variables), then for each
field $F$ the set
\begin{equation}\tag{1}\label{1}
\bigl\{c\in F^M : f_0(c,X)\in\bigl(f_1(c,X),\dots,f_n(c,X)\bigr)F[X]\bigr\}
\end{equation}
is a {\it constructible}\/ subset of $F^M$, i.e., a boolean combination of
algebraic subsets of $F^M$. 
Theorem~C below can be seen as an analogue for polynomials with
coefficients in $\Z$.
Results on dependence on parameters such as this are most conveniently
(and accurately) expressed using the terminology of mathematical logic:
for example, Hermann's Theorem asserts that 
the set \eqref{1} above is definable by a quantifier-free formula in
the language ${\cal L}=\{{0},{1},{+},{-},{\cdot}\}$ of rings, for all
fields $F$. (See \cite{Chang-Keisler} or \cite{Hodges} 
for the basic notions of first-order logic and model theory.)

Before we can state the next theorem, we have to introduce some more notations.
If $a,b\in\Z$ are not both zero, we let $(a:b):=a/\gcd(a,b)$, where 
$\gcd(a,b)$ is a positive generator of the ideal $(a,b)$ of $\Z$.
We also put $(0:0):=1$. Moreover, we define
a relation $\rad$ on pairs $(a,b)$ of integers as follows:
$\rad(a,b)$ holds if and only if $b$ divides $a^n$, for some $n\in\N:=
\{0,1,2,\dots\}$.
Let ${\cal L}_{\rad}$ be the expansion of the language 
$\cal{L}$ by a binary function symbol $(\ :\ )$ and a binary
predicate symbol $\rad$. We construe the ring $\Z$ as
${\cal L}_{\rad}$-structure by interpreting the ring symbols as usual and
the symbols $(\ :\ )$ and $\rad$ as described above.

\begin{theoremC} 
Let $f_0(C,X),f_1(C,X),\dots,f_n(C,X)\in\Z[C,X]$.
The set
\begin{equation}\label{2}\tag{2}
\bigl\{c\in\Z^M : f_0(c,X)\in\bigl(f_1(c,X),\dots,f_n(c,X)\bigr)\Z[X]\bigr\}
\end{equation}
is definable by a quantifier-free formula in the language ${\cal L}_{\rad}$.
\end{theoremC}

It follows that for fixed $f_0(C,X),\dots,f_n(C,X)\in\Z[C,X]$, one can
decide in polynomial time whether a tuple $c\in\Z^M$ is in the set \eqref{2}.
(This is well-known for $N=0$, see, e.g., \cite{vdDries-Holly}.)
The quantifier-free formula in question can even be constructed 
from the $f_j$'s by a primitive recursive algorithm.

Here is an analogue of Theorem~C for homogeneous equations.
We say that a term $\tau(C,X)$ in a language $\cal{L}^*$ extending the language
${\cal L}=\{{0},{1},{+},{-},{\cdot}\}$ of rings
is {\bf polynomial in $X$} if $\tau(C,X)=
f\bigl(\tau^*(C),X\bigr)$ for some polynomial $f\in\Z[V,X]$, where 
$V=(V_1,\dots,V_L)$ is a tuple of new variables
and $\tau^*(C)$ an $L$-tuple of $\cal{L}^*$-terms in the 
variables $C$. 
(In other words, the extra function symbols in ${\cal L}^*\setminus
\cal{L}$ are applied
only to subterms of $\tau$ not involving the $X$-variables.)
We let ${\cal L}_{\gcd}$ be the sublanguage $\{ {0},{1},{+},{\cdot}, { (\ :\ ) }\}$ of ${\cal L}_{\rad}$.

\begin{theoremD} 
Let $f_1(C,X),\dots,f_n(C,X)\in\Z[C,X]$. There exists a
finite family $\bigl\{\varphi^{(\lambda)}(C)\bigr\}_{\lambda\in\Lambda}$ of 
quantifier-free ${\cal L}_{\gcd}$-formulas
and for each $\lambda\in\Lambda$ finitely many $n$-by-$1$ 
column vectors
$$y^{(\lambda, 1)}(C,X),\dots,y^{(\lambda, K)}(C,X) \qquad (K\in\N)$$
whose entries are ${\cal L}_{\gcd}$-terms, polynomial
in $X$, such that for all $c\in\Z^M$ we have
$\Z\models\bigvee_{\lambda\in\Lambda} \varphi^{(\lambda)}(c)$,
and if $\lambda\in\Lambda$ is such that
$\Z\models\varphi^{(\lambda)}(c)$, then
$$y^{(\lambda, 1)}(c,X),\dots,y^{(\lambda, K)}(c,X)\in\Z[X]^n$$
generate the $\Z[X]$-module of solutions in $\Z[X]$ to the homogeneous equation
$$f_1(c,X)y_1+\cdots+f_n(c,X)y_n=0.$$
\end{theoremD}

In fact, the $\varphi^{(\lambda)}$ and the $y^{(\lambda, k)}$ 
can be effectively constructed from $f_1,\dots,f_n$.
Theorems C and D (suitably adapted)
remain true in the more general setting of B\'ezout domains.
It should be remarked that in contrast to Theorem~D, it is not possible 
in general to obtain a parametric solution $(y_1,\dots,y_n)\in \Z[X]^n$
(uniform in the parameters $c\in\Z^M$) to an inhomogeneous linear equation 
\begin{equation}\label{3}\tag{3}
f_0(c,X)=f_1(c,X)y_1+\cdots+f_n(c,X)y_n,
\end{equation}
even for the case $f_0=1$. 
More precisely, by Theorem~B (or the example in Section~6 of
\cite{maschenb-ideal2})
there do not exist finitely many
$n$-tuples $\bigl(\tau_{1k}(C,X),\dots,\tau_{nk}(C,X)\bigr)$ 
of terms in a language $\cal{L}^*\supseteq \cal{L}$
such that $\Z$ can be expanded to an ${\cal L}^*$-structure, 
each $\tau_{ik}(C,X)$ is polynomial in $X$, and
such that if $c\in\Z^M$ with
$$1\in\bigl(f_1(c,X),\dots,f_n(c,X)\bigr)\Z[X],$$ then
$\bigl(\tau_{1k}(c,X),\dots,\tau_{nk}(c,X)\bigr)\in\Z[X]^n$ is a solution
to \eqref{3}, for some $k$. 

\subsection*{Organization of the paper}
Sections~\ref{Preliminaries-Section} and \ref{Coherent-Section}
contain preliminary material. Beside fixing notations, 
we introduce a tool from first-order logic, namely
the persistence of Horn formulas under direct products (or more generally:
reduced products). This allows us to shorten some arguments in later sections
(although it is not strictly speaking necessary).
In Section~\ref{Coherent-Section} we discuss coherent modules and rings.
Most of the material is standard, but we emphasize issues of
uniformity and definability. In Section~\ref{Homogeneous-Section} we
study bounds for homogeneous systems of linear equations. We introduce
a notion (super coherence) related to the notion of
``stable coherence'' from \cite{Glaz} and prove Theorems~A
and D. In Section~\ref{Inhomogeneous-Section} we prove Theorems~B and
C. The theorems in Sections~\ref{Homogeneous-Section} and 
\ref{Inhomogeneous-Section} can be employed to obtain uniformity and
definability results for various properties of ideals
and algebraic constructions
in polynomial rings. In Section~\ref{Prime-Section} we illustrate this by
means of defining the primeness of an ideal. 
In an appendix (Section~\ref{Appendix})
we give yet another application of the material in 
Section~\ref{Inhomogeneous-Section} and obtain a characterization of
Jacobson domains among Noetherian domains inspired by a characterization
of Noetherian domains with the ``Skolem property'' in \cite{Frisch}.

\subsection*{Acknowledgments} Parts of this paper derive from
my Ph.D. thesis \cite{maschenb-thesis}. I would like
to thank Lou van den Dries for his guidance during writing of that thesis,
in particular for pointing out Proposition~\ref{NilProp} below.  
I am also grateful
to the University of California at Berkeley and the Mathematical Sciences
Research Institute, where this article was written.

\section{Preliminaries}\label{Preliminaries-Section}

In this section we collect some definitions and notations used in the sequel. 
The reader may glance over this part and come back to it for
reference when necessary. We also recall some basic facts about Horn
formulas which will be handy in Sections~\ref{Coherent-Section} and
\ref{Homogeneous-Section}.

\subsection*{Rings, ideals and modules}
Let $R$ be a ring (throughout: commutative with a unit element $1$). 
We write $(r_1,\dots,r_n)R$ for the ideal generated in $R$ by elements
$r_1,\dots,r_n$; we omit $R$ if it is clear from the context.
The localization $S^{-1}R$, where $S$ is the set of non-zero-divisors
of $R$, is called the ring of fractions of $R$, denoted by $\Frac(R)$.
For submodules $M$, $M'$  of an $R$-module we define the ideal
$$M': M:=\bigl\{a\in R:am\in M' \text{ for all $m\in M$}\bigr\},$$
of $R$.
Given an $R$-module $M$ and $(f_1,\dots,f_n)\in M^n$, 
the set of solutions
in $R^n$ to the homogeneous system of linear equations 
$y_1f_1+\cdots+y_nf_n=0$ 
is an $R$-submodule of $R^n$, which we call 
the (first) {\bf module of syzygies of $(f_1,\dots,f_n)$}. 
If $M=R^m$ and $f_1,\dots,f_n\in R^m$ are the column vectors of
a matrix $A\in R^{m\times n}$, we denote
the module of syzygies of $(f_1,\dots,f_n)$ by $\Sol_R(A)$
(the module of solutions to the system of homogeneous linear equations
$Ay=0$). If
$I$ is an ideal of $R$, then 
$$\sqrt{I}=\bigl\{r\in R : \text{$r^n\in I$ for some $n>0$}\bigr\}$$
is the nilradical of $I$. 
%and
%$$\Jac(I)=\bigl\{r\in R : \text{$1+rs$ is a unit mod $I$, for all $s\in R$}\bigr\}$$
%the Jacobson radical of $I$. Equivalently, $\sqrt{I}$ is the intersection
%of all prime ideals of $R$ containing $I$, and $\Jac(I)$ the intersection of all
%maximal ideals of $R$ containing $I$. 
We let $\Nil(R):=\sqrt{(0)}$, the nilradical of $R$.
%and 
%$\rad(R):=\Jac\bigl((0)\bigr)$, the
%Jacobson radical of $R$.

%If
%$R$ is coherent (e.g., if $R$ is Noetherian), 
%then $\Sol_R(A)$ is finitely generated.
%For submodules $M$, $M'$  of an $R$-module we write
%$$(M': M):=\bigl\{a\in R:am\in M' \text{ for all $m\in M$}\bigr\},$$
%an ideal of $R$ (containing the annihilator of $M$).

\subsection*{Polynomials}
Unless otherwise noted,
by $X=(X_1,\dots,X_N)$ we always denote a tuple of $N$ distinct 
indeterminates, where $N\in\N$. The (total)
degree of a polynomial $0\neq f\in R[X]=R[X_1,\dots,X_N]$ is denoted
by $\deg(f)$.
By convention $\deg(0) := -\infty$ where $-\infty<\N$.
We extend this notation to finite tuples
$f=(f_1,\dots,f_n)$ of polynomials in $R[X]$ by setting
$\deg(f):=\max_j \deg(f_j)$ (the degree of $f$).

\subsection*{Semihereditary rings}
A ring $R$ is called {\bf hereditary} if every ideal of $R$ is projective,
and {\bf semihereditary} if every finitely generated ideal of $R$ is
projective. 
If $R$ is a domain, then $R$ is hereditary if and only if $R$ is a 
Dedekind domain, and $R$ is semihereditary if and only if $R$ is 
a Pr\"ufer domain. (\cite{Glaz}, p.~27.) 
Every semihereditary ring is reduced, and every von Neumann regular ring 
is semihereditary. 
A ring $R$ is semihereditary if and only if
$\Frac(R)$ is von Neumann regular and $R_{\frak m}$ is a valuation ring
for every maximal ideal $\frak m$ of $R$. (\cite{Glaz}, Corollary~4.2.19.)
If $R$ is hereditary, then $R_{\frak p}$
is a discrete valuation ring (DVR) for every prime ideal $\frak p$ of $R$.

\subsection*{Reduced products}
Let $\cal F$ be a filter on $\N$, i.e., a collection of non-empty subsets
of $\N$ closed under taking finite intersections and supersets.
For every $k\in\N$ let $R^{(k)}$ be a ring. 
The {\bf reduced product} 
$R^*=\prod_{k\in\N} R^{(k)}/\cal F$ of 
the family $\{R^{(k)}\}_{k\in\N}$ over $\cal F$ is the 
ring $R/I$, where $I$ is the ideal
of $R=\prod_{k\in\N} R^{(k)}$ 
consisting of all sequences $a=(a^{(k)})\in R$ with 
$\{k:a^{(k)}=0\}\in\cal F$. We write $a\mapsto a/{\cal F}:=a+I$ for
the canonical homomorphism $R\to R^*=R/I$, and extend it in the usual
manner to a homomorphism
$R^n\to (R^*)^n$ (for $n\in\N$) which we denote in the same way.
If ${\cal F}=\{\N\}$, then $a\mapsto a/\cal F$ is an
isomorphism $R\to R^*$.
Now suppose in addition 
that for every $k\in\N$ we are given an $R^{(k)}$-module
$M^{(k)}$. Similarly as above, we then define the reduced product 
$M^*=\prod_{k\in\N} M^{(k)}/\cal F$ of
$\{M^{(k)}\}_{k\in\N}$ over the filter $\cal F$ by $M^*=M/N$, where
$N$ is the submodule of the $R$-module 
$M=\prod_k M^{(k)}$ consisting of all sequences
$m=(m^{(k)})\in M$ with $\{k:m^{(k)}=0\}\in\cal F$. Then $M^*$ is an
$R^*$-module.

\subsection*{Horn formulas}
Let $R$ be a ring and $M$ an $R$-module. We construe $M$ as a two-sorted
structure (in the sense of model theory) in the following way: The two
sorts are the {\bf ring sort} with underlying set $R$ and variables
$r,s,\dots$, and the {\bf group sort} with underlying sort $M$ and
variables $x,y,\dots$. The corresponding two-sorted language
${\cal L}_{\operatorname{mod}}$ of modules is the disjoint union of:
\begin{enumerate}
\item the language ${\cal L}=\{{0},{1},{+},{-},{\cdot}\}$ of
rings,  interpreted in the obvious way in $R$;
\item the language $\{{0},{+},{-}\}$ of additive groups,
interpreted in the obvious way in $M$;
\item a binary function symbol $\cdot$,  interpreted as
scalar multiplication $(r,x)\mapsto r\cdot x\colon R\times M\to M$.
\end{enumerate}
A {\bf basic Horn formula} is an ${\cal L}_{\operatorname{mod}}$-formula of the form 
$$\sigma_1= 0\ \&\ \cdots\ \&\ \sigma_p= 0 \ \rightarrow\ 
\tau_1= 0\ \&\ \cdots\ \&\ \tau_q= 0$$ 
where $p$ and $q$ are natural numbers, $q\geq 1$, 
and $\sigma_i=\sigma_i(r,x)$, $\tau_j=\tau_j(r,x)$ are 
${\cal L}_{\bmod}$-terms in two collections of distinct indeterminates
$r=(r_1,\dots,r_m)$ (ranging over $R$) and 
$x=(x_1,\dots,x_n)$ (ranging over $M$).
We allow the case $p=0$, in which case the formula in
question is just $\tau_1=\cdots=\tau_q= 0$.
A {\bf Horn formula} is an ${\cal L}_{\operatorname{mod}}$-formula 
consisting of a finite (possibly empty)
string of quantifiers, followed by a conjunction of basic Horn formulas.
A Horn formula that is an ${\cal L}_{\operatorname{mod}}$-sentence
is called a {\bf Horn sentence.}

For each $k\in\N$ let $M^{(k)}$ be a module over the ring $R^{(k)}$.
Let $\cal F$ be a filter
on $\N$, and let $R$, $M$, $R^*$ and $M^*$ be as above.
For every $k$ let 
$a^{(k)}=\bigl(a_1^{(k)},\dots,a_m^{(k)}\bigr)\in (R^{(k)})^m$ and
$b^{(k)}=\bigl(b_1^{(k)},\dots,b_n^{(k)}\bigr)\in (M^{(k)})^n$. We put
$a_i=(a_i^{(k)})\in R$, $b_i=(b_i^{(k)})\in M$ and
$a=(a_1,\dots,a_m)$, $b=(b_1,\dots,b_n)$.
The following is a special case of a
fundamental theorem about Horn formulas due to Chang.
(In the case of a direct product, i.e., $\cal{F}=\{\N\}$,
it was first proved by Horn.)

\begin{theorem}\label{Horn}
For any Horn formula $\varphi(r,x)$ as above,
$$\bigl\{k\in\N : M^{(k)}\models\varphi(a^{(k)},b^{(k)})\bigr\}\in\cal F
\quad\Rightarrow\quad
M^*\models\varphi(a/{\cal F},b/{\cal F}).$$
\end{theorem}

We omit the straightforward proof of this theorem
(see, e.g., \cite{Hodges}, Theorem~9.4.3) and instead,
as an illustration for its usefulness, apply it to
reprove a well-known algebraic fact:

\begin{lemma}\label{Product}
Every reduced product of a family of semihereditary rings is semihereditary.
\end{lemma}
\begin{proof}
A ring $R$ is semihereditary if and only if for all $n\geq 1$ and all
$f_1,\dots,f_n\in R$ the following holds, with 
$I:=(f_1,\dots,f_n)R$ and $\phi\colon R^n\to I$, 
$\phi(a_1,\dots,a_n)=a_1f_1+\cdots+a_nf_n$:
There exist $n^2$ elements $y_{ij}\in R$ such that
the map $\psi\colon I\to R^n$ given by $\psi(f_i)= 
(y_{i1},\dots,y_{in})$ is well-defined and satisfies 
$\phi\circ\psi=\id_I$.
For given $n$, this statement can be easily 
formalized as a Horn sentence.
The claim now follows from Theorem~\ref{Horn} (in the case where
$M^{(k)}=R^{(k)}$ for all $k$).
\end{proof}

Theorem~\ref{Horn} also admits a converse: for any sentence $\psi$
in the  language of modules which is preserved under reduced 
products of families of
modules there is a Horn sentence which is equivalent to $\psi$, in any
module. This much deeper fact,
due to Galvin and Keisler, will not be used here; 
see \cite{Chang-Keisler}, Theorem~6.2.5.

\section{Coherent Modules and Coherent Rings}\label{Coherent-Section}

In this section, $R$ always denotes 
a ring. An $R$-module $M$ is {\bf finitely presented}
(sometimes also called {\bf finitely related}) if there exists an
exact sequence $F_1 \to F_0 \to M\to 0$ of $R$-linear maps,
where $F_0,F_1$ are finitely
generated free $R$-modules. 
A finitely generated $R$-module $M$ is called {\bf coherent}
if every finitely generated submodule of $M$ is finitely presented.
Every finitely generated submodule of a coherent module is itself a
coherent module. If $R$ is Noetherian, then every finitely generated
$R$-module is coherent.

We call a finitely generated $R$-module
$M$ {\bf $\alpha$-uniformly coherent}, where $\alpha\colon\N\to\N$ is a
function, if for every $n\in\N$
the kernel of every $R$-module homomorphism $R^n\to M$ is generated by at 
most $\alpha(n)$ many elements. (Equivalently, the syzygies of every
element of $M^n$ can be generated by $\alpha(n)$ elements of $R^n$, for
all $n\in\N$.) In this case, 
we call the function $\alpha$ a 
{\bf uniformity function}\/ for $M$. We
say that $M$ is {\bf uniformly coherent} if it is $\alpha$-uniformly
coherent for some uniformity function $\alpha$;
clearly uniformly coherent $\Rightarrow$ coherent. (Uniformly coherent modules
were first defined and studied by Soublin \cite{Soublin-1}; 
see also \cite{Glaz}, \cite{Glaz-History}.)

We say that an $R$-module $M$ is {\bf $m$-generated} (for $m\in\N$) if
it is generated by $m$ elements.
Being $m$-generated and $\alpha$-uniformly coherent 
is a property of a module (for given $m$ and given uniformity function 
$\alpha$) which is preserved under taking reduced products. More precisely:

\begin{prop}\label{Prod}
Let $m\in\N$ and $\alpha\colon\N\to\N$ be a function. Let
$\{R^{(k)}\}_{k\in\N}$ be a family of rings,
and for every $k\in\N$ let $M^{(k)}$ be an $m$-generated and
$\alpha$-uniformly coherent $R^{(k)}$-module.
Then for every filter $\cal F$  on $\N$, $\prod_k M^{(k)}/\cal F$ 
is a module over $\prod_k R^{(k)}/\cal F$ which is $m$-generated and 
$\alpha$-uniformly coherent.
\end{prop}
\begin{proof}
%For each $k\in\N$ let $g_1^{(k)},\dots,g_m^{(k)}\in M^{(k)}$ be generators for the 
%$R^{(k)}$-module $M^{(k)}$. Then $(g_1^{(k)}),\dots,(g_m^{(k)})$ generate the
%module $\prod_k M^{(k)}$ over $\prod_k R^{(k)}$. Hence the elements
%$g_j=\bigl(g_j^{(k)}\bigr)/\cal F$, $j=1,\dots,m$, generate the $R^*$-module $M^*$. 
%Let $n\in\N$, $n\geq 1$, and $f_1,\dots,f_n\in M^*$ be given. 
%Write $f_i=\bigl(f_i^{(k)}\bigr)/\cal F$ with 
%$f_i^{(k)}\in M^{(k)}$. Then for every $k\in\N$,
%there exist $$y^{(j,k)}=\left(y_1^{(j,k)},\dots,y_n^{(j,k)}\right)\in (R^{(k)})^n,
%\qquad j=1,\dots,\alpha(n),$$ which generate
%the submodule of $(R^{(k)})^n$ consisting of the solutions to the homogeneous
%linear equation $$y_1f_1^{(k)}+\cdots+y_nf_n^{(k)}=0.$$ It is straightforward 
%to check that then $$y^{(j)}=\left(\bigl(y_1^{(j,k)}\bigr)/{\cal F},\dots,
%\bigl(y_n^{(j,k)}\bigr)/{\cal F}\right)\in (R^*)^n,\qquad
%j=1,\dots,\alpha(n),$$ generate the $R^*$-module of
%solutions to the homogeneous linear equation
%$$y_1f_1+\cdots+y_nf_n=0.$$ This shows that $M^*$ is $\alpha$-coherent
%as required.
By Theorem~\ref{Horn}, since
the condition that a given module is
$m$-generated and $\alpha$-uniformly 
coherent can be expressed by a Horn sentence.
\end{proof}

\begin{cor}\label{Uniform-Coherence}
For an $R$-module $M$, an integer $m\geq 0$ 
and a function $\alpha\colon\N\to\N$,
the following are equivalent:
\begin{enumerate}
\item $M$ is $m$-generated and $\alpha$-uniformly coherent.
\item For every filter $\cal F$ on $\N$, 
$M^\N/\cal F$ is an $m$-generated and
$\alpha$-uniformly coherent $R^\N/\cal F$-module.
%\item For some ultrafilter $\cal U$ on $\N$, 
%$M^\N/\cal U$ is an $\alpha$-uniformly coherent  $R^\N/\cal U$-module. \qed
\item $M^\N$ is an $m$-generated and $\alpha$-uniformly
coherent $R^\N$-module.
\end{enumerate}
\end{cor}

The proposition above also yields a 
characterization of uniform coherence due to Soublin:

%\begin{proof}
%The implication (1)~$\Rightarrow$~(2) follows from 
%the proposition, 
%and (2)~$\Rightarrow$~(3) by taking ${\cal F}=\{\N\}$ in (2). 
%Suppose $M^\N$ is $m$-generated and $\alpha$-coherent. 
%Then $M$ is also $m$-generated.
%We write $\delta_M$ and $\delta_R$
%for the diagonal embeddings $M\to M^\N$ and $R\to R^\N$, respectively.
%Let $f_1,\dots,f_n\in M$. Then the syzygies of 
%$\bigl(\delta_M(f_1),\dots,\delta_M(f_n)\bigr)\in (M^\N)^n$ are generated by $\alpha(n)$
%elements.
%Any syzygy $(y_1,\dots,y_n)$ of $(f_1,\dots,f_n)$ in $R^n$ gives rise to
%a syzygy $\bigl(\delta_R(y_1),\dots,\delta_R(y_n)\bigr)$ of
%$\bigl(\delta_M(f_1),\dots,\delta_M(f_n)\bigr)\in (M^\N)^n$
%in $(R^\N)^n$. It follows that $M$ is $\alpha$-coherent.
%\end{proof}

\begin{cor}\label{Uniform-Coherence-2}
The following are equivalent, for an $R$-module $M$:
\begin{enumerate}
\item $M$ is finitely generated and uniformly coherent.
\item For every filter $\cal F$ on $\N$,
$M^\N/\cal F$ is a finitely generated coherent $R^\N/\cal F$-module.
\item $M^\N$ is a finitely generated coherent $R^\N$-module.
\end{enumerate}
\end{cor}
\begin{proof}
The implication (1)~$\Rightarrow$~(2) follows from the proposition, 
and (2)~$\Rightarrow$~(3) by taking ${\cal F}=\{\N\}$ in (2). 
It remains to show
(3)~$\Rightarrow$~(1). So assume that
$M^\N$ is an $m$-generated coherent module over $R^\N$, for
some $m\in\N$. Then $M$ is an $m$-generated $R$-module. 
Suppose for a contradiction that
$M$ is not uniformly coherent, that is,
there is an integer $n\in\N$ with the following property: 
For every $k\in\N$ there is $\bigl(f_1^{(k)},\dots,f_n^{(k)}\bigr)\in M^n$
whose syzygies cannot be generated by $k$ elements. Let
$f_i=(f_i^{(k)}) \in M^\N$ for $i=1,\dots,n$. Then the $R^\N$-module 
of syzygies of $(f_1,\dots,f_n)\in (M^\N)^n$ is not finitely generated,
contradicting the coherence of $M^\N$.
\end{proof}

A ring $R$ is called {\bf coherent} if it is coherent as a module over
itself, that is, if every finitely generated ideal of $R$ is finitely
presented. The following characterizations of coherence are due to Chase
\cite{Chase}; for a proof see \cite{Glaz}, p.~45--47.

\begin{theorem}\label{Theorem-Chase}
The following are equivalent, for a ring $R$:
\begin{enumerate}
\item $R$ is a coherent ring.
\item Every finitely presented $R$-module is coherent.
\item Every direct product of flat $R$-modules is flat.
\item For every non-empty set $\Lambda$, the $R$-module $R^\Lambda$ is flat.
\item For every finitely generated ideal $I$ of $R$ and every $a\in R$,
the ideal $I:(a)$ is finitely generated.
\item For every $a\in R$, the ideal $(0):(a)$ of $R$ is finitely generated, and
the intersection of two finitely generated ideals of $R$ is finitely
generated ideal.
\end{enumerate}
\end{theorem}

An ideal $I$ of $R$ is called {\bf nilpotent} if there exists
an integer $m\geq 1$ such that $I^m=\{0\}$, and the smallest such $m$
is called the {\bf index of nilpotency} of $I$. 
Here are some sufficient
conditions which ensure the preservation of 
coherence under ring extensions and quotients. (See \cite{Glaz}, 
Theorem~4.1.1.)

\begin{prop}\label{Descent}
Let $\phi\colon R\to S$ be a ring homomorphism making $S$ into a finitely
presented $R$-module. 
\begin{enumerate}
\item If $R$ is a coherent ring, then $S$ is a coherent ring.
\item If $\ker\phi$ is a finitely presented nilpotent ideal of $R$ and
$S$ is a coherent ring, then $R$ is a coherent ring.
\end{enumerate}
\end{prop}

A ring $R$ is called {\bf $\alpha$-uniformly coherent} 
if it is $\alpha$-uniformly coherent as a module over itself, and
{\bf uniformly coherent} if it is $\alpha$-uniformly coherent for
some $\alpha\colon\N\to\N$. By Corollary~\ref{Uniform-Coherence-2}, $R$ is
uniformly coherent if and only if $R^\N$ is coherent. 
Noetherian rings are rarely uniformly coherent: A Noetherian ring $R$
is $\alpha$-uniformly coherent if and only if $\dim R\leq 2$ and
$R_{\frak m}$ is $\alpha$-uniformly coherent,
for every maximal ideal $\frak m$ of $R$. In this case one can take
$\alpha(n)=n+2$. (See \cite{Glaz}, Corollary~6.1.21.) 

The condition of $\alpha$-coherence only concerns syzygies of tuples of
elements of $R$. However, 
it implies the existence of a finite bound on the number
of generators for the syzygies of tuples of elements of $R^m$ for $m>1$:

\begin{lemma}\label{Coherence-Matrix-Lemma}
For all integers $m,n>0$ and all $m\times n$-matrices $A$
with entries in an $\alpha$-uniformly coherent ring $R$, the module
$\Sol_R(A)$
of solutions to the homogeneous system of linear equations
$Ay=0$ is generated by $\alpha^m(n)$ solutions.
\end{lemma}
\begin{proof}
We proceed by induction on $m$, the case $m=1$ just being the definition
of $\alpha$-coherence. Suppose $m>1$, and
let $n$ be  positive integer, $R$ $\alpha$-uniformly coherent, 
and $A$ an $m\times n$-matrix
with entries from $R$. Let $(f_1,\dots,f_n)$ be the first
row of $A$ and $A'$ be the matrix consisting of the last $m-1$ rows of $A$.
Let $z_1,\dots,z_\alpha$  be generators for the
syzygies of $(f_1,\dots,f_n)$, where $\alpha=\alpha(n)$.  
Consider the $z_i$ as column vectors and
let $B=A'\cdot (z_1,\dots,z_\alpha)$, an $(m-1)\times \alpha$-matrix 
with entries in $R$. 
The solutions to $Ay=0$ are in one-to-one correspondence with the
solutions to $Bu=0$: every solution $u=(u_1,\dots,u_\alpha)^{\tr}
\in R^\alpha$ to $Bu=0$ gives rise to a solution $y=(y_1,\dots,y_n)^{\tr}
\in R^n$ 
to $Ay=0$ by setting $y=\sum_i u_iz_i$, and every solution to $Ay=0$ arises
in this way. By inductive hypothesis, there are 
$\alpha^{m}(n)$ generators for the module of solutions
to $Bu=0$, giving rise to as many
generators for the module of solutions to $Ay=0$.
\end{proof}

\begin{definition}
Let $\cal C$ be a class of rings. We say that $\cal C$ is {\bf 
$\alpha$-uniformly coherent} if every member of $\cal C$ is $\alpha$-uniformly
coherent. We call $\cal C$ {\bf uniformly coherent} if it is $\alpha$-uniformly
coherent for some uniformity function $\alpha$.
\end{definition}

The following lemma gives a criterion for a class of rings to be
uniformly coherent. 

\begin{lemma}\label{Direct-Product}
Suppose that $\cal C$ is a class of rings which is closed under direct 
products. Then $\cal C$ is uniformly coherent if and only if every 
$R\in\cal C$ is coherent.
\end{lemma}
\begin{proof}
The ``only if'' direction is trivial. The proof of the ``if'' direction
is similar to the proof of (3)~$\Rightarrow$~(1) in 
Corollary~\ref{Uniform-Coherence-2}:
Suppose for a contradiction
that every $R\in\cal C$ is coherent, but $\cal C$ is not
uniformly coherent. Then there exists some $n\in\N$ such that for
every $k\in\N$ there is an $R^{(k)}\in\cal C$ and $\bigl(f_1^{(k)},\dots,
f_n^{(k)}\bigr)\in (R^{(k)})^n$ whose syzygies in $(R^{(k)})^n$ cannot be
generated by $k$ elements. Let $R^*=\prod_k R^{(k)}$ and
$f_i=(f_i^{(k)})\in R^*$. Then $R^*\in\cal C$, so $R^*$ is coherent. But
the module of syzygies of $(f_1,\dots,f_n)$ in $(R^*)^n$ is not
finitely generated, a contradiction.
\end{proof}

The typical example of a uniformly coherent class of rings is
the class of semihereditary rings:

\begin{lemma}\label{Semi-uc}
Every semihereditary ring is uniformly coherent
with uniformity function $\alpha(n)=n$.
\end{lemma}
\begin{proof}
Let $R$ be a semihereditary ring and let $f_1,\dots,f_n\in R$. 
We have to show that the syzygies of
$(f_1,\dots,f_n)$ are generated by $n$ elements of $R^n$. The finitely 
generated ideal $I:=(f_1,\dots,f_n)$ of $R$ is projective. So the short
exact sequence $$0\to K:=\ker\phi\to R^n\overset{\phi}{\to} I\to 0,$$ where
$\phi(a_1,\dots,a_n)=a_1f_1+\cdots+a_nf_n$ for $(a_1,\dots,a_n)\in R^n$,
splits. Hence $K$ is a direct summand of $R^n$, and thus
generated by $n$ elements.
\end{proof}

Other examples for uniformly coherent classes of rings can be
obtained from rings of finite rank: Inspired by a definition of 
I.~S.~Cohen \cite{Cohen} we say that
a ring $R$ {\bf has finite rank} if for some natural number $k>0$,
every finitely generated ideal of $R$ is generated by $k$ elements. 
(In \cite{Cohen} this definition is only made for Noetherian $R$.)
We call the smallest integer $k>0$ with this property the 
{\bf rank} of $R$.
Equivalently, a ring $R$ has rank $k$ if every ideal of
$R$ which is generated by $k+1$ elements is generated by $k$ elements, but
there exists a finitely generated 
ideal of $R$ which cannot be generated by fewer than $k$ elements. 
For example, the
domains of rank $1$ are exactly the B\'ezout domains.
Any reduced product of
a family of rings of rank $\leq k$ has itself rank $\leq k$, by 
Theorem~\ref{Horn}.

\begin{lemma}\label{Finite-Rank}
Let $R$ be a coherent ring of rank $k$. Then every finitely generated
submodule of $R^n$ can be generated by $nk$ elements.
\end{lemma}
%\begin{proof}
%Let $M$ be a finitely generated submodule of $R^n$, say $M=\operatorname{im}\phi$ with
%an $R$-linear map $\phi\colon R^m\to R^n$ given by $\phi(x)=Ax$, where
%$A$ is an $n\times m$-matrix with entries from 
%$R$. For each $i=1,\dots,n$ consider the ideal
%$$I_i := \bigl\{ r_i : \text{$r_ie_i+\cdots+r_ne_n\in M$
%for some $r_{i+1},\dots,r_n\in R$}\bigr\}$$
%of $R$. (Here $e_1,\dots,e_n$ denote the standard basis
%vectors of $R^n$.) Let $a_1,\dots,a_n$ denote the row vectors of $A$,
%and let $A^{(i)}$ be the $(i-1)\times n$-matrix with row vectors
%$a_1,\dots,a_{i-1}$. Since $R$ is coherent, there
%exist finitely many column vectors $b_1,\dots,b_L\in R^m$ which generate
%the module of solutions to the linear equation $A^{(i)}x=0$ in $R^m$.
%The elements $a_ib_l$ ($l=1,\dots,L$) generate the ideal $I_i$.
%Since $R$ has rank $k$, there exist $r_{i1},\dots,r_{ik}\in R$ with
%$I_i=(r_{i1},\dots,r_{ik})R$.
%Choose $v_{ij}\in M$ and $s_{ij}\in R$ with 
%$$v_{ij}=r_{ij}e_i+s_{i+1,j}e_{i+1}+\cdots+s_{n,j}e_n.$$
%It is easy to check that the $v_{ij}$ ($1\leq i\leq n$,
%$1\leq j\leq k$) generate $M$.
%\end{proof}

This lemma appears in \cite{Cohen}, for Noetherian $R$; the proof given
there goes through for $R$ coherent.

\begin{cor}\label{Finite-Rank-Cor}
The class of coherent rings of rank $k$
is uniformly coherent with uniformity function $\alpha(n)=nk$.
\end{cor}
\begin{proof}
Let $R$ be a coherent ring of rank $k$, and let
$\phi\colon R^n\to R$ be an $R$-linear map. 
Since $R$ is coherent, $\ker\phi$ is
finitely generated. By the previous lemma, $\ker\phi$
can be generated by $nk$ elements.
\end{proof}

Let us mention some classes of coherent rings with finite rank.
First note that
every Artinian ring $R$ has finite rank, equal to the length of $R$.
A Noetherian domain 
$R$ has finite rank if and only if $\dim R\leq 1$. 
(\cite{Cohen}, Theorem~9.) See \cite{Bass} for information about
Noetherian domains of rank $2$.

\begin{prop}
Let $R$ be a ring of finite Krull dimension $d$.
\begin{enumerate}
\item If each localization of $R$ has rank $\leq k$, then $R$ has
rank $\leq d+k$.
\item If each localization of $R$ is uniformly coherent with common
uniformity function $n\mapsto \alpha(n)$, then $R$ is uniformly coherent with
uniformity function $n\mapsto d+\alpha(n)$.
\end{enumerate}
\end{prop}
\begin{proof}
We use the following fact, which is the
culmination of work of Forster \cite{Forster}, Swan \cite{Swan-1}, Eisenbud and
Evans \cite{Eisenbud-Evans} (for Noetherian rings) and Heitmann 
\cite{Heitmann-1}, \cite{Heitmann-2} (in the general case):
Let $M$ be a 
finitely generated $R$-module, and $k\in\N$. 
If for each prime ideal $\frak p$ of $R$,
the $R_{\frak p}$-module $M_{\frak p}=M\otimes_R R_{\frak p}$ can be
generated by $k$ elements, then the $R$-module
$M$ can be generated by $d+k$ elements.

The first part of the proposition now 
follows immediately. For the second
part, let $\phi\colon R^n\to R$ be an $R$-linear map. Then
$\ker (\phi\otimes_R\id_{R_{\frak p}})=(\ker \phi)\otimes_R R_{\frak p}$ 
can be generated by $\alpha(n)$ 
elements, for every prime ideal $\frak p$ of $R$. Hence
$\ker \phi$ can be generated by $d+\alpha(n)$  elements.
\end{proof}

\begin{remarks}
By part (1) in the proposition
it follows that a hereditary ring $R$ has rank $\leq 2$, since $\dim R=1$ and 
each localization of $R$, being a DVR, has rank $\leq 1$.
In contrast, there exist Pr\"ufer domains of finite rank $>2$ 
\cite{Schuelting}, and even of infinite rank \cite{Swan-2}. 
As to (2), note that the assumption $\dim R<\infty$
cannot be dropped: there exists a ring $R$ all of whose localizations are
valuation rings, but $R$ is not
coherent (\cite{Glaz}, p.~54).
\end{remarks}

\subsection*{Model-theoretic aspects}
Let ${\cal L}^*$ be a language extending the language 
${\cal L}=\{0,1,+,-,\cdot\}$ of rings, and let $\cal C$ be a class of 
${\cal L}^*$-structures whose $\cal L$-reducts are rings.
Fix $\alpha\colon\N\to\N$, and
suppose that for every integer $n\geq 1$ there is a finite family
$\bigl\{\varphi^{(\lambda)}(C)\bigr\}_{\lambda\in\Lambda}$
of ${\cal L}^*$-formulas $\varphi^{(\lambda)}(C)$,
where $C=(C_1,\dots,C_n)$ is an $n$-tuple of distinct variables, and for
each $\lambda\in\Lambda$ finitely many $n$-by-$1$ column vectors
$$y^{(\lambda,1)}(C),\dots,y^{(\lambda,\alpha(n))}(C)$$
whose entries are ${\cal L}^*$-terms, with the following properties:
For every $R\in\cal C$ and every
$f=(f_1,\dots,f_n)\in R^n$, we have
\begin{enumerate}
\item $R\models\bigvee_{\lambda\in\Lambda}\varphi^{(\lambda)}(f)$;
\item if $\lambda\in\Lambda$
is such that $R\models\varphi^{(\lambda)}(f)$, then the vectors
$$y^{(\lambda,1)}(f),\dots,y^{(\lambda,\alpha(n))}(f)\in R^n$$
generate the $R$-module of syzygies of $f$. 
\end{enumerate}
In particular, $R$ is $\alpha$-uniformly coherent.
By the proof of Lemma~\ref{Coherence-Matrix-Lemma}, 
for all integers $m,n>0$ and every
$m\times n$-matrix  $A=(a_{ij})\in R^{m\times n}$ there
exists a similar parametrization of generators for the $R$-module
$\Sol_R(A)$ by $\alpha^m(n)$ many column vectors whose entries are
${\cal L}^*$-terms, which is
uniform in $R$ and the $a_{ij}$.
%finite family 
%$\bigl\{\varphi^{(\gamma)}_{m,n}(C)\bigr\}_{\gamma\in\Gamma}$
%of ${\cal L}^*$-formulas $\varphi^{(\gamma)}_{m,n}(C)$,
%where here $C=(C_{ij})$ is an $mn$-tuple of distinct variables, and for
%each $\gamma\in\Gamma$ finitely many column vectors
%$$y^{(\gamma,1)}(C),\dots,y^{(\gamma,\alpha^m(n))}(C)$$
%whose entries are ${\cal L}^*$-terms, 
%such that:
%\begin{enumerate} 
%\item For every $R\in\cal C$ and every $m\times n$-matrix $A=(a_{ij})\in 
%R^{m\times n}$, we have
%$R\models\bigvee_\lambda\psi^{(\lambda)}(A)$;
%\item if $\lambda\in\Lambda$
%is such that $R\models\psi^{(\lambda)}(A)$, then 
%the $R$-module $\Sol_R(A)$ is generated by
%the vectors
%$$z^{(\lambda,1)}(A),\dots,z^{(\lambda,\alpha^m(n))}(A)\in R^n.$$
%\end{enumerate}
Moreover, if the $\varphi^{(\lambda,j)}$ can be chosen quantifier-free 
(for all $n$), then the corresponding formulas describing the parametrization
of the generators for $\Sol_R(A)$ can also be chosen
quantifier-free.

We now consider an important example.
Let ${\cal L}_{\gcd}=\bigl\{{0},{1},{+},{-},{\cdot},{(\ :\ )}\bigr\}$
be the language
obtained by augmenting the language ${\cal L}$ of rings by a binary 
function symbol $(\ :\ )$. Every B\'ezout domain $R$ can be construed as an
${\cal L}_{\gcd}$-structure by interpreting $(\ :\ )$ in the
following way. For $a,b\in R$ let
$\gcd(a,b)$ be a generator of the ideal
$(a,b)$ generated by $a$ and $b$,
chosen in such a way that 
\begin{equation}\label{gcd-0}
\gcd(a,b)=\gcd(b,a). 
\end{equation}
For example, if $R$ is a valuation
ring, then we may define
$\gcd(a,b):=b$ if $b\neq 0$ divides $a$ and $\gcd(a,b):=a$ otherwise;
for $R=\Z$ we may choose $\gcd(a,b)$ to be the unique non-negative
generator of $(a,b)$.
The element
$$(a:b) := \begin{cases}
\frac{a}{\gcd(a,b)} &\text{if $a\neq 0$ or $b\neq 0$} \\
1                   &\text{otherwise}                  
\end{cases}$$
generates the ideal
$$(a):(b)=\bigl\{c\in R:bc\in (a)\bigr\}.$$
Note that by \eqref{gcd-0}, for all $a,b\in R$:
\begin{equation}\label{gcd-1}
b\cdot(a:b)=a\cdot(b:a)
\end{equation}
and hence
\begin{equation}\label{gcd-2}
(a:b)\cdot (bc:ac) = (b:a)\cdot (ac:bc) \quad\text{for all non-zero $c\in R$.}
\end{equation}
%Inductively, we also define
%$$\gcd(a_1,\dots,a_n) := \gcd\bigl(\gcd(a_1,\dots,a_{n-1}),a_n\bigr)\qquad
%\text{for $a_1,\dots,a_n\in R$, $n> 2$,}$$
%a generator for the ideal $(a_1,\dots,a_n)$.
%
Let now $R$ be a ring, $a=(a_1,\dots,a_n)\in R^n$, and
consider the  homogeneous linear equation
\begin{equation}\label{BezHom}
a_1y_1+\cdots+a_ny_n=0.
\end{equation}
We have:

\begin{lemma}\label{BezLemma}
If $R$ is a B\'ezout domain and $a\neq 0$, then the
module of solutions in $R^n$ to \eqref{BezHom}
is generated by the special solutions
$$y^{(i,j)}=\bigl[0,\dots,0,(a_j:a_i),0,\dots,0,-(a_i:a_j),0,\dots,0\bigr]^{\operatorname{tr}}\in R^n, \quad
1\leq i<j\leq n.$$
\end{lemma}

This is a consequence of the following general observation (valid for
any ring $R$):

\begin{prop}
Suppose that
$\lambda_1,\dots,\lambda_n\in R$ are such that
$$u=\lambda_1 a_1+\cdots+\lambda_n a_n$$
is a unit in $R$.
Then the module of solutions in $R^n$ to \eqref{BezHom}
is generated by the $n$ special solutions
$$y^{(i)}=\left[\lambda_1 a_i,\dots,\lambda_{i-1}a_i,-\sum_{k\neq i}\lambda_k
a_k,\lambda_{i+1}a_i,\dots,\lambda_n a_i\right]^{\operatorname{tr}},
\quad i=1,\dots,n.$$
\end{prop}
\begin{proof}
We have
$$ay^{(i)}=\sum_{j\neq i} (\lambda_ja_i)a_j-\left(\sum_{k\neq i}\lambda_k
a_k\right)a_i=0,$$
so $y^{(i)}$ is a solution to \eqref{BezHom}, for $i=1,\dots,n$. Let $y=(y_1,\dots,y_n)^{\operatorname{tr}}\in R^n$ be any solution to \eqref{BezHom}.
The $i$th component of the vector 
$y_1y^{(1)}+\cdots+y_ny^{(n)}$ is given by
$$\sum_{j\neq i} y_j\cdot (\lambda_i a_j)-y_i\cdot\left(\sum_{k\neq i}\lambda_k a_k\right)=
-\lambda_i a_i y_i-y_i(1-\lambda_ia_i)=-uy_i.$$
Hence $$y=-u^{-1}\bigl(y_1y^{(1)}+\cdots+y_ny^{(n)}\bigr),$$ showing that the $y^{(j)}$ generate
the module of solutions to \eqref{BezHom} in $R^n$.
\end{proof}

We now prove Lemma~\ref{BezLemma}; so 
suppose that $R$ is a B\'ezout domain. Clearly the $y^{(i,j)}$ are
solutions to \eqref{BezHom}, by \eqref{gcd-1}.
Let $0\neq d\in R$ be a generator for $(a_1,\dots,a_n)R$. 
Then the linear homogeneous equation
$$\frac{a_1}{d}y_1+\cdots+\frac{a_n}{d}y_n=0$$
over $R$ has the same solutions in $R^n$ as \eqref{BezHom}, and 
for all $1\leq i,j\leq n$
there exists a unit $u$ of $R$ such that
$(a_i/d:a_j/d)=u\cdot (a_i:a_j)$ and
$(a_j/d:a_i/d)=u\cdot (a_j:a_i)$, by \eqref{gcd-2}. So replacing
$a_i$ by $a_i/d$ for all $i$, if necessary, we may assume that
$$1=\lambda_1 a_1+\cdots+\lambda_n a_n \qquad\text{for some $\lambda_1,\dots,
\lambda_n\in R$.}$$
Let $y^{(1)},\dots,y^{(n)}$ be as in the proposition. One shows easily that
$$y^{(i)}=\sum_{j=1}^{i-1} \lambda_i\gcd(a_i,a_j)\cdot y^{(j,i)} -
      \sum_{j=i+1}^n \lambda_j\gcd(a_i,a_j)\cdot y^{(i,j)}$$
for all $i=1,\dots,n$. Therefore, since the $y^{(1)},\dots,y^{(n)}$ 
generate the
solution module of \eqref{BezHom}, so do the $y^{(i,j)}$ ($1\leq i<j\leq n$). \qed

\medskip

Lemma~\ref{BezLemma} and the discussion
at the beginning of this subsection (applied to ${\cal L}^*={\cal L}_{\gcd}$)
yield the following fact.  Here $C=(C_1,\dots,C_M)$.

\begin{cor}\label{BezCor}
Let $A(C)\in\Z[C]^{m\times n}$.
One can construct elementary recursively \rom{(}from $A$\rom{)} a finite
family $\bigl\{\varphi^{(\lambda)}(C)\bigr\}_{\lambda\in\Lambda}$ of
quantifier-free ${\cal L}_{\gcd}$-formulas $\varphi^{(\lambda)}(C)$ and for
each $\lambda\in\Lambda$ finitely many $n$-by-$1$ column vectors
$$y^{(\lambda,1)}(C),\dots,y^{(\lambda,n)}(C)$$
whose entries are ${\cal L}_{\gcd}$-terms, such that for all B\'ezout domains 
$R$ and
$c\in R^M$, we have $R\models\bigvee_\lambda \varphi^{(\lambda)}(c)$, and
if $\lambda\in\Lambda$ is such that $R\models\varphi^{(\lambda)}(c)$, then
the vectors 
$$y^{(\lambda,1)}(c),\dots,y^{(\lambda,n)}(c)\in R^n$$
generate the $R$-module $\Sol_{R}\bigl(A(c)\bigr)$. 
\end{cor}

This corollary slightly improves
\cite{vdDries-Holly}, Corollary~5.4, where a similar parametrization
was given using terms in a larger language.

\subsection*{A flatness result}
We finish this section by proving a fact about subrings of
direct products of rings (Corollary~\ref{Almost-All-Zero})
which will be used in the next section.

Let $\{R^{(k)}\}_{k\in\N}$ be a family of rings and
$R^*=\prod_k R^{(k)}$ its direct product.
Given $k\in\N$ we identify $r\in R^{(k)}$ with the sequence $(r^{(l)})\in R^*$
given by $r^{(l)}=0$ for $l\neq k$ and $r^{(l)}=r$ for $l=k$. In this way,
$R^{(k)}$ becomes an ideal of $R^*$. 
%We consider the direct sum $\bigoplus_k R^{(k)}$ 
%of $\{R^{(k)}\}$ as an ideal 
%of the direct product $R^*=\prod_k R^{(k)}$
%of $\{R^{(k)}\}$ in the usual way.
We write $\pi^{(k)}\colon R^*\to R^{(k)}$
for the projection onto the 
$k$-th component: $\pi^{(k)}(y)=y^{(k)}$ for $y=(y^{(k)})\in R^*$. 
We extend $\pi^{(k)}$ in the natural way 
to a ring homomorphism $(R^*)^n\to (R^{(k)})^n$, denoted by 
the same symbol.
Let $S$ be a subring of $R^*$.

\begin{lemma}
Let $M$ be a finitely generated $S$-submodule of $S^n$ and
$M^*$ be an $R^*$-submodule of $(R^*)^n$ with $M^*\supseteq M$. If
$\pi^{(k)}(M^*) = \pi^{(k)}(M)$ for all $k$, then $M^*=R^*M$ 
\rom{(}$=\text{the $R^*$-submodule of $M^*$ generated by $M$}$\rom{)}.
\end{lemma}
\begin{proof}
Let $g_1,\dots,g_m\in M$ be generators for the $S$-module $M$, and let $y\in M^*$. 
Then for every $k\in\N$ we can write 
$y^{(k)}=a_1^{(k)}g_1^{(k)}+\cdots+a_m^{(k)}g_m^{(k)}$ 
for some $a_i^{(k)}\in R^{(k)}$.
Putting $a_i=(a_i^{(k)})\in R^*$ we obtain $y=a_1g_1+\cdots+a_mg_m\in R^*M$ as
required.
\end{proof}

\begin{cor}\label{Almost-All-Zero}
Suppose that $S\supseteq \bigoplus_{k} R^{(k)}$ is coherent.
Then $R^*$ is a flat $S$-module.
\end{cor}
\begin{proof}
Let $M^*$ be the module of syzygies in $(R^*)^n$ of a tuple
$(f_1,\dots,f_n)\in S^n$, so
$M=M^*\cap S^n$ is a finitely generated $S$-module.
We have to show $R^*M=M^*$, and hence,
by the lemma above, that 
$\pi^{(k)}(M)=\pi^{(k)}(M^*)$ for every $k$.
For this, let $y=(y_1,\dots,y_n)\in M^*$, so $f_1y_1+\cdots+f_ny_n=0$.
Let $y^{(k)}=(y_1^{(k)},\dots,y_n^{(k)})\in (R^{(k)})^n$. Then
$f_1y_1^{(k)}+\cdots+f_ny_n^{(k)}=0$ in $R^*$ and hence
$y^{(k)}\in M^*\cap (R^{(k)})^n$. Since $R^{(k)}\subseteq S$
this yields $y^{(k)}\in M$ and hence
$\pi^{(k)}(y)=y^{(k)}=\pi^{(k)}(y^{(k)})\in \pi^{(k)}(M)$ as required.
\end{proof}

\section{Homogeneous Linear Equations in Polynomial Rings}\label{Homogeneous-Section}

In this section we will be concerned with the existence of uniform bounds for
the degrees of generators for syzygy modules over polynomial rings. We
define a {\it super coherent}\/ class of rings to be one 
for which such bounds exist. 
This notion is related to ``stable coherence'' introduced in \cite{Glaz}. 
We show that the class of semihereditary rings is super coherent, yielding
Theorem~A from the Introduction. We also prove Theorem~D and discuss some
strengthenings of Theorem~A.

\subsection*{Stable coherence and super coherence}
A ring $R$ is called {\bf stably coherent} if for
every $N\geq 0$ the ring of polynomials $R[X_1,\dots,X_N]$ over $R$
is coherent. We say that a class $\cal C$ of rings is stably coherent
if every $R\in\cal C$ is stably coherent. 
For example, the class of Noetherian rings
is stably coherent, by virtue of the Hilbert Basis Theorem.
There exist coherent rings $R$ which are not stably coherent
\cite{Soublin-2}.
We have the following theorem proved by Vasconcelos,
after a conjecture by Sabbagh (\cite{Sabbagh}, p.~502). For an efficient proof
based on work of Alfonsi see \cite{Glaz}, Chapter~7.

\begin{theorem}\label{Vasconcelos}
The class of semihereditary rings is stably coherent.
\end{theorem}

For the purpose of this section 
we introduce a notion related to stable coherence:

\begin{definition}\label{Definition-Supercoherence}
Let $\alpha\colon\N\to\N$ and $\beta\colon\N^3\to\N$.
We call a ring $R$ {\bf $(\alpha,\beta)$-super coherent} if 
$R$ is $\alpha$-uniformly coherent, and for given $N,d,n\in\N$ and
$f_1,\dots,f_n\in R[X]=R[X_1,\dots,X_N]$ of degree at most $d$, 
every solution to the homogeneous linear equation
$$f_1y_1+\cdots+f_ny_n=0$$
is a linear combination of solutions of degree at most $\beta(N,d,n)$.
We say that $R$ is {\bf super coherent} if it is $(\alpha,\beta)$-super
coherent for some functions $\alpha$ and $\beta$ as above.
\end{definition}

\begin{remarks}
Let $R$ be $(\alpha,\beta)$-super coherent. 
The localization $R_U$ of $R$ at a multiplicative subset $U$ of $R$
is $(\alpha,\beta)$-super coherent. 
If $I$ is a finitely generated ideal of $R$, then $R/I$ 
is super coherent. If $R$ is a faithfully flat extension of an
$\alpha$-uniformly coherent subring $S$, then $S$ is $(\alpha,\beta)$-super
coherent. 
(These facts are
immediate consequences of the definition.)
\end{remarks}

A super coherent ring is stably coherent. In fact:

\begin{lemma}\label{SupercoherentIsStablycoherent}
Given $\alpha\colon\N\to\N$ and $\beta\colon\N^3\to\N$ there exists a 
function $\gamma\colon\N^3\to\N$ with the
property that for all  $(\alpha,\beta)$-super coherent rings $R$ and
all $f_1,\dots,f_n\in R[X]$ 
of degree $\leq d$, the module of solutions in $R[X]$ 
to the homogeneous linear equation
\begin{equation}\label{H1}
f_1y_1+\cdots+f_ny_n=0
\end{equation}
is generated by $\gamma(N,d,n)$ many solutions of degree 
$\leq\beta(N,d,n)$.
\end{lemma}
\begin{proof}
The $R$-module of solutions to
the homogeneous linear equation \eqref{H1} in $R[X]^n$ 
which have degree at most
$\beta=\beta(N,d,n)$ is isomorphic to the module of solutions to a 
certain system of
$m':=\binom{N+\beta+d}{N}$ homogeneous linear equations over $R$ in
$n':=n\cdot\binom{N+\beta}{N}$ indeterminates. Hence by 
Lemma~\ref{Coherence-Matrix-Lemma} the former module
can be generated by $\gamma(N,d,n):=
\alpha^{m'}(n')$ many elements.
These elements will then also generate the $R[X]$-module of solutions
to \eqref{H1} in $R[X]^n$.
\end{proof} 

If $R$ is $(\alpha,\beta)$-super coherent and of finite rank $k$, then 
we can take $\alpha(n)=nk$ (Corollary~\ref{Finite-Rank-Cor}), and
the function
$(N,d,n)\mapsto \beta(N,d,n)$ and hence also the function $(N,d,n)\mapsto
\gamma(N,d,n)$ can be chosen so as to not 
depend on $n$. For if we have a bound
$\beta=\beta(N,d,n)$ 
for $n=\binom{N+d}{N}\cdot k$, then this $\beta$ will also be a bound for
all other values of $n$. (By Lemma~\ref{Finite-Rank}.)

Lemma~\ref{SupercoherentIsStablycoherent} extends to systems of
homogeneous linear equations:

\begin{cor}
For any given $N,d,m,n\in\N$, $m,n\geq 1$
there exist natural numbers $\beta_m=\beta_m(N,d,n)$ and $\gamma_m=
\gamma_m(N,d,n)$ with the
property that for all $(\alpha,\beta)$-super coherent rings $R$ 
and all $m\times n$-matrices $A$
with entries in $R[X]=R[X_1,\dots,X_N]$ of degree $\leq d$, the module
$\Sol_{R[X]}(A)$ 
is generated by $\gamma_m$ elements of degree $\leq\beta_m$.
\end{cor}
\begin{proof}
We proceed by induction on $m$, similar to the proof of 
Lemma~\ref{Coherence-Matrix-Lemma}.
Clearly $\beta_1=\beta$ and $\gamma_1=\gamma$ as in 
Lemma~\ref{SupercoherentIsStablycoherent} work for $m=1$. Suppose $m>1$, and
let $N,d,n\in\N$ with $n\geq 1$. Let $R$ be $(\alpha,\beta)$-super coherent
and $A$ an $m\times n$-matrix
with entries from $R[X]=R[X_1,\dots,X_N]$. Let $(f_1,\dots,f_n)$ be the first
row of $A$ and $A'$ be the matrix consisting of the last $m-1$ rows of $A$.
Let the column vectors $z_1,\dots,z_\gamma$  of degree $\leq \beta$ 
generate the
syzygies of $(f_1,\dots,f_n)$, and put $B=A'\cdot(z_1,\dots,z_\gamma)$, an $(m-1)\times \gamma$-
matrix with entries in $R[X]$. Here $\beta=\beta(N,d,n)$, and
$\gamma=\gamma(N,d,n)$ is as in Lemma~\ref{SupercoherentIsStablycoherent}.
Every solution $u=(u_1,\dots,u_\gamma)^{\tr}
\in R[X]^\gamma$ to $Bu=0$ gives rise to a solution 
$y=(y_1,\dots,y_n)^{\tr}\in R[X]^n$ 
to $Ay=0$ by setting $y=\sum_i u_iz_i$. This yields a
one-to-one correspondence between the solutions to $Bu=0$ and
the solutions to $Ay=0$. The degrees of the entries of $B$ are bounded 
from above by $\beta+d$.
By inductive hypothesis, there are $\gamma_{m-1}(N,\beta+d,\gamma)$ generators
of degree $\leq \beta_{m-1}(N,\beta+d,\gamma)$ for the module of solutions
to $Bu=0$. These give rise to  $\gamma_{m-1}(N,\beta+d,\gamma)$
generators of degree $\leq \beta\cdot  \beta_{m-1}(N,\beta+d,\gamma)$
for the module of solutions to $Ay=0$.
Hence we can take $\beta_m(N,d,n)=\beta\cdot  \beta_{m-1}(N,\beta+d,\gamma)$ 
and $\gamma_m(N,d,n)=\gamma_{m-1}(N,\beta+d,\gamma)$.
\end{proof}
\begin{remark}
The proof shows that if $(N,d,n)\mapsto\beta(N,d,n)$ 
does not depend on $n$, then 
$\beta_m(N,d,n)$ and $\gamma_m(N,d,n)$ 
can be chosen independent of $n$, for all $m\geq 1$.
\end{remark}

Let $R$ be a ring. We say that an $R[X]$-submodule of
$R[X]^m$ is {\bf of type $(n,d)$} (where $n,d\in\N$) if it is generated by $n$
elements of degree $\leq d$. The previous corollary and
standard arguments (see, e.g., 
\cite{maschenb-ideal2}, proof of Proposition~4.7) yield:

\begin{cor}\label{Type-Cor}
Given $\alpha\colon\N\to\N$ and $\beta\colon\N^3\to\N$ there exists a 
function $\tau\colon\N^4\to\N^2$ with the following properties:
if $R$ is an $(\alpha,\beta)$-super coherent ring and
$M,M'$ are finitely generated submodules of the free
$R[X]$-module $R[X]^m$ of type $(n,d)$, then the $R[X]$-module
$M\cap M'$ and the ideal $M':M$ of $R[X]=R[X_1,\dots,X_N]$
are of type $\tau(N,d,m,n)$. If $\beta(N,d,n)$ does not depend on $n$, then
$\tau(N,d,m,n)$ also does not depend on $n$. 
\end{cor}

A class $\cal C$ of rings is called {\bf $(\alpha,\beta)$-super coherent} 
if every ring $R\in\cal C$ is $(\alpha,\beta)$-super coherent. We say that
$\cal C$ is {\bf super coherent} if $\cal C$ is $(\alpha,\beta)$-super coherent
for some $\alpha,\beta$ as above.
The main result of this section is the following:

\begin{theorem}\label{SuperIsStable}
Let $\cal C$ be a class of rings which is closed under direct products. Then
$\cal C$ is super coherent if and only if $\cal C$ is stably coherent.
\end{theorem}

In particular, it then follows that a ring $R$ is super coherent if and only
if $R^\N$ is stably coherent.
The theorem together with
Lemma~\ref{Product} and Theorem~\ref{Vasconcelos} implies:

\begin{cor}\label{SemiIsSuper}
The class of semihereditary rings is super coherent. 
\end{cor}

\begin{remarks}
The previous corollary implies Theorem~A stated in the Introduction.
By Lemma~\ref{Semi-uc} above and Theorem~4.1 of \cite{maschenb-ideal2},
the class of hereditary rings is $(\alpha,\beta)$-super coherent with
$\alpha(n)=n$ and $\beta(N,d)=(2d)^{2^{O(N^2)}}$.
\end{remarks}

Before we give a proof of Theorem~\ref{SuperIsStable}, we
establish some auxiliary facts. 
Let $\cal F$ be a filter on $\N$ and $\{R^{(k)}\}_{k\in\N}$ a family of
rings indexed by $\N$. Let $R=\prod_k R^{(k)}$ and 
$R^* = \prod_k R^{(k)}/\cal F$.
The following observation (not used later on) might give an
indication why direct products rather than reduced products (or ultraproducts,
as in \cite{vdDries-Schmidt}, say) play
the most prominent role in our investigations:

\begin{lemma}\label{Flatness}
The canonical homomorphism $a\mapsto a/{\cal F}\colon R\to R^*$ is flat.
%\rom{(In particular, if $R$ is $(\alpha,\beta)$-super coherent, so is $R^*$.)}
\end{lemma}
\begin{proof}
We have $R^* = R/I$, where $I$ denotes the ideal of $R$ consisting of all
sequences $a=(a^{(k)})\in R$ with $\{k:a^{(k)}=0\}\in{\cal F}$.
By \cite{Glaz}, Theorem~1.2.15, $R/I$ is a flat $R$-module if and only 
if for every $a\in I$ there exists a $c\in I$ with $(1-c)a=0$.
To see this, let $a=(a^{(k)})\in I$, so
$\Delta:=\{k:a^{(k)}=0\}\in{\cal F}$. Define
$c=(c^{(k)})\in R$ by $c^{(k)}=0$ if $k\in\Delta$ and $c^{(k)}=1$ if
$k\notin\Delta$. Then $c\in I$ and $(1-c)a=0$ as required. 
\end{proof}

Let $R^*[X]$ be the ring of polynomials in indeterminates 
$X=(X_1,\dots,X_N)$ with coefficients from $R^*$, and put
$R[X]^*= \prod_k R^{(k)}[X]/\cal F$. We have a natural
embedding of $R^*$-algebras $R^*[X]\to R[X]^*$
induced by $X_i\mapsto X_i/{\cal F}\in R[X]^*$ for $i=1,\dots,N$. 
We consider $R^*[X]$ as a subring of $R[X]^*$ via this
embedding. Note that if $\cal F=\{\N\}$, then 
$R^*[X]=\prod_k R^{(k)}[X]$ becomes identified
in this way with the $R$-subalgebra of the direct product
$R[X]^*=\prod_k R^{(k)}[X]$ consisting of all
sequences $(f^{(k)})$ of polynomials whose degrees are bounded,
that is, such that there exists $d\in\N$ with
$\deg f^{(k)}\leq d$ for all $k$.

\begin{lemma}\label{SuperCoh}
The following are equivalent, for a uniformly coherent class $\cal C$ of rings:
\begin{enumerate}
\item $\cal C$ is super coherent.
\item For every family $\{R^{(k)}\}_{k\in\N}$ of rings in $\cal C$, every
filter $\cal F$ on $\N$ and
every $N\in\N$, $R[X]^*$ is flat over
$R^*[X]$, where $X=(X_1,\dots,X_N)$.
\item For every family $\{R^{(k)}\}_{k\in\N}$ of rings in $\cal C$ and every
$N\in\N$, the ring
$\prod_k R^{(k)}[X]$ is flat over $\left(\prod_k R^{(k)}\right)[X]$, with $X=(X_1,\dots,X_N)$.
\end{enumerate}
\rom{(In particular, a uniformly coherent ring $R$ is super coherent 
if and only if $R[X]^\N$ is flat over $R^\N[X]$.)}
\end{lemma}
\begin{proof}
Suppose $\cal C$ is $(\alpha,\beta)$-super coherent, let
$N,d,n\in\N$, $n\geq 1$ be fixed, and let $\gamma=\gamma(N,d,n)$ be as in 
Lemma~\ref{SupercoherentIsStablycoherent} above. 
Let $f_1(C,X),\dots,f_n(C,X)\in\Z[C,X]$ be
general polynomials of degree $d$, where $X=(X_1,\dots,X_N)$, and 
$C=(C_1,\dots,C_M)$ are
parametric variables. It is easy to write down a Horn formula $\varphi(C)$
(where $C$ is considered as a tuple of variables of the ring sort)
which, for a given ring $R$ and $c\in R^M$, 
holds in $R[X]$ (considered as a module over
itself) for $c$ exactly if there exist $\gamma$
solutions to the equation 
$$f_1(c,X)y_1+\dots+f_n(c,X)y_n=0$$ in $R[X]$ of degree
$\leq\beta$ from which every solution to this homogeneous
equation in $R[X]$ can be obtained as an $R[X]$-linear combination. 
Hence (2) is a consequence of Theorem~\ref{Horn}.
The implication
(2)~$\Rightarrow$~(3) follows by taking ${\cal F}=\{\N\}$.
For (3)~$\Rightarrow$~(1)
suppose for a contradiction that (3) holds but $\cal C$ 
is not super coherent. So there exist 
$N,d,n\in\N$ with $n\geq 1$, 
and for every $k\in\N$ a ring $R^{(k)}\in\cal C$ and polynomials
$f_1^{(k)},\dots,f_n^{(k)}\in R^{(k)}[X]=R^{(k)}[X_1,\dots,X_N]$ of degree at most $d$
such that the module of solutions in $R^{(k)}[X]$ to the homogeneous
linear equation 
\begin{equation}\label{Equ-k}
f_1^{(k)}y_1+\cdots+f_n^{(k)}y_n=0
\end{equation} 
cannot be generated by
elements of degree $\leq k$, that is, there exists a column vector
$$y^{(k)}=\left[y_1^{(k)},\dots,y_n^{(k)}\right]^{\tr}\in R^{(k)}[X]^n$$ with
$$f_1^{(k)}y_1^{(k)}+\cdots+f_n^{(k)}y_n^{(k)}=0$$ which is not
an $R^{(k)}[X]$-linear combination of solutions of degree $\leq k$.
Put $$R^*:=\prod_k R^{(k)}, \qquad R[X]^* := \prod_k R^{(k)}[X].$$
Write each polynomial $f_i^{(k)}$ as $$f_i^{(k)}=\sum_\nu a_{i,\nu}^{(k)}X^\nu,$$
where the sum ranges over all
$\nu=(\nu_1,\dots,\nu_N)\in\N^N$ and 
$a_{i,\nu}^{(k)}\in R^{(k)}$, $X^\nu=X_1^{\nu_1}\cdots X_N^{\nu_N}$.
We have $a_{i,\nu}^{(k)}=0$ if
$\abs{\nu}:=\nu_1+\cdots+\nu_N>d$. Hence
$$f_i = \sum_{\nu} a_{i,\nu}X^\nu \in R^*[X],$$ where
$a_{i,\nu} = (a_{i,\nu}^{(k)})_{k\in\N}\in R^*$, with
$a_{i,\nu} = 0$ if $\abs{\nu}>d$.
The column vector
$y=[y_1,\dots,y_n]^{\tr}$, where $y_i=(y_i^{(k)})\in R[X]^*$, is a solution 
to the homogeneous linear equation
$$f_1y_1+\cdots+f_ny_n=0.$$
Since $R[X]^*$ is flat over $R^*[X]$, $y$ is an $R[X]^*$-linear 
combination of certain solutions $z_1,\dots,z_m\in\bigl(R^*[X]\bigr)^n$.
Let $k$ be an integer larger than the degrees of $z_1,\dots,z_m$. It
follows that $y^{(k)}=\pi^{(k)}(y)$ is a linear combination of the solutions 
$$\pi^{(k)}(z_1),\dots,\pi^{(k)}(z_m)\in \bigl(R^{(k)}[X]\bigr)^n$$ 
to \eqref{Equ-k} which have degree $\leq k$. This is a
contradiction to the choice of $y^{(k)}$. 
\end{proof}

We now prove Theorem~\ref{SuperIsStable}. Let $\cal C$ be a class of rings
which is closed under direct products.
We have already remarked that any super coherent ring is stably coherent.
Suppose conversely that $\cal C$ is stably coherent.
By Lemma~\ref{Direct-Product}, $\cal C$ is uniformly coherent.
In order to show that $\cal C$ is super coherent, we have to prove
that for every family $\{R^{(k)}\}_{k\in\N}$ of rings in $\cal C$ and every
integer $N\geq 0$,
$R[X]^*=\prod_k R^{(k)}[X]$ is flat over $R^*[X]=
\left(\prod_k R^{(k)}\right)[X]$, with $X=(X_1,\dots,X_N)$.
(By Lemma~\ref{SuperCoh}.)
Clearly the subring $R^*[X]$ of $R[X]^*$ contains 
$\bigoplus_k \bigl(R^{(k)}[X]\bigr)$.
Since $\cal C$ is closed under direct products
and stably coherent, $R^*[X]$ is coherent.
The claim now follows from 
Corollary~\ref{Almost-All-Zero}. \qed

\begin{remark}
The hypothesis of Theorem~\ref{SuperIsStable} is cannot be dropped, as
the class ${\cal G}_2$ of coherent rings of global dimension $2$ 
shows. (See, e.g., \cite{Glaz} for the definition of the global dimension of a ring.)
By a theorem of Greenberg and Vasconcelos \cite{Greenberg-Vasconcelos},
${\cal G}_2$ is stably coherent.
However, there exist rings in ${\cal G}_2$ which are not super coherent:
For example let $R=\Q[[U,V]]$, where $U$, $V$ are distinct indeterminates,
and consider the ideals
$$I= (U-VX), \qquad J_d = (UV^d,U^d-2V^d)$$
of the polynomial ring $R[X]$,
where $X$ is a single indeterminate and $d\geq 4$. 
Then $I\cap J_d$ cannot be generated by polynomials of degree $<d$,
see \cite{Soublin-2}. By Corollary~\ref{Type-Cor} it follows that $R$ is not
super coherent.
\end{remark}

\subsection*{Model-theoretic consequences}
Let now
$$A(C,X)=
\bigl(a_{ij}(C,X)\bigr)_{\substack{1\leq i\leq m\\
1\leq j\leq n}}$$ be an $m\times n$-matrix with entries $a_{ij}(C,X)
\in \Z[C,X]$, where $C=(C_1,\dots,C_M)$.
The following fact is an immediate consequence of Corollaries~\ref{BezCor}
and \ref{SemiIsSuper}.

\begin{cor}\label{DefBezout-Cor}
There exists a finite
family $\bigl\{\varphi^{(\lambda)}(C)\bigr\}_{\lambda\in\Lambda}$ of
quantifier-free ${\cal L}_{\gcd}$-formulas $\varphi^{(\lambda)}(C)$ and for
each $\lambda\in\Lambda$ finitely many $n$-by-$1$ column vectors
$$y^{(\lambda,1)}(C,X),\dots,y^{(\lambda,K)}(C,X)\qquad
(K\in\N)$$
whose entries are ${\cal L}_{\gcd}$-terms in the variables $(C,X)$,
polynomial in $X$, such that for all B\'ezout domains 
$R$ and
$c\in R^M$, we have $R\models\bigvee_\lambda \varphi^{(\lambda)}(c)$, and
if $\lambda\in\Lambda$ is such that $R\models\varphi^{(\lambda)}(c)$, then
the vectors 
$$y^{(\lambda,1)}(c,X),\dots,y^{(\lambda,K)}(c,X)\in R[X]^n$$
generate the $R[X]$-module $\Sol_{R[X]}\bigl(A(c,X)\bigr)$. 
\end{cor}

\begin{remark}
The case $m=1$, $R=\Z$ of the corollary yields Theorem~D in the
Introduction.
\end{remark}

Let ${\cal L}_{\div}$ be ${\cal L}=\{{0},{1},{+},{-},{\cdot}\}$
augmented by a binary predicate $|$,
to be interpreted in every ring $R$ as divisibility, that is:
$a|b :\Longleftrightarrow ac=b \text{ for some $c\in R$.}$
For the special case of valuation rings, Corollary~\ref{DefBezout-Cor} gives:

\begin{cor}\label{DefBezout-Cor-2}
There exists a finite
family $\bigl\{\psi^{(\lambda)}(C)\bigr\}_{\lambda\in\Lambda}$ of
quantifier-free ${\cal L}_{\div}$-formulas $\psi^{(\lambda)}(C)$ and for
each $\lambda\in\Lambda$ a polynomial $\delta^{(\lambda)}(C)\in\Z[C]$ and
finitely many column vectors
$$z^{(\lambda,1)}(C,X),\dots,z^{(\lambda,K)}(C,X)\qquad
(K\in\N)$$
whose entries $z_j^{(\lambda)}(C,X)$ are polynomials in $\Z[C,X]$, 
such that for all
valuation rings $R$ and $c\in R^M$,
we have $R\models\bigvee_\lambda \psi^{(\lambda)}(c)$, and
if $\lambda\in\Lambda$ is such that $R\models\psi^{(\lambda)}(c)$, then
$\delta^{(\lambda)}(c)\neq 0$ divides 
$z_j^{(\lambda,k)}$ in $R[X]$, for every $j,k$, and
$$\bigl(z^{(\lambda,1)}/\delta^{(\lambda)}\bigr)(c,X),\dots,\bigl(z^{(\lambda,K)}/\delta^{(\lambda)}\bigr)(c,X)\in R[X]^n$$
generate  $\Sol_{R[X]}\bigl(A(c,X)\bigr)$. 
\end{cor}

\begin{remark}
The remark following Corollary~\ref{SemiIsSuper} shows that from $A(C,X)$ one
can  (elementary recursively) {\it construct}\/ a finite family 
$\bigl\{\varphi^{(\lambda)}(C)\bigr\}$ of
quantifier-free ${\cal L}_{\gcd}$-formulas and corresponding
column vectors $y^{(\lambda,k)}$ which satisfy the property expressed in
Corollary~\ref{DefBezout-Cor} {\it for every principal ideal domain
\rom{(}PID\rom{)} $R$.}\/
Similarly, from $A(C,X)$ one can explicitly construct the objects
$\psi^{(\lambda)}$, $\delta^{(\lambda)}$ and $z^{(\lambda,k)}$ 
having the properties stated in the previous corollary {\it for every
DVR $R$.}
\end{remark}

Let $A'(C,X)\in\Z[C,X]^{m\times n'}$ where $n'\in\N$, $n'\geq 1$. 
For each ring $R$ and $c\in R^M$,
we may consider  the $R[X]$-submodules $M(c,X)$ and $M'(c,X)$ of the
free $R[X]$-module $R[X]^m$ generated by the columns of  $A(c,X)$ and
$A'(c,X)$, respectively.
Corollary~\ref{DefBezout-Cor-2} immediately implies the 
uniformity of certain module-theoretic operations on $M(c,X)$ and $M'(c,X)$:

\begin{cor}\label{DefBezout-Cor-3}
There exists a finite
family $\bigl\{\theta^{(\lambda)}(C)\bigr\}_{\lambda\in\Lambda}$
of quantifier-free ${\cal L}_{\div}$-formulas, and for each
$\lambda\in\Lambda$ tuples a polynomial $\delta^{(\lambda)}(C)\in\Z[C]$,
an $m\times K$-matrix 
$B^{(\lambda)}(C,X)\in\Z[C,X]^{m\times K}$
and polynomials
$$u^{(\lambda,1)}(C,X),\dots,u^{(\lambda,K)}(C,X),$$
for some integer $K\geq 1$,
with the following property: for every valuation ring $R$ and $c\in R^M$, 
we have
$R\models\bigvee_{\lambda\in\Lambda} \theta^{(\lambda)}$, and if
$R\models\theta^{(\lambda)}(c)$,
then $\delta^{(\lambda)}(c)\neq 0$ divides all entries of $B^{(\lambda)}(c,X)$ and
the $u^{(\lambda,k)}(c,X)$ in $R[X]$, 
the $R[X]$-module $$M(c,X)\cap M'(c,X)$$ is generated by the columns of 
the matrix $B^{(\lambda)}(c,X)/\delta^{(\lambda)}(c)$, and
the ideal $$M'(c,X) : M(c,X)$$ of $R[X]$ is generated by 
$$u^{(\lambda,1)}(c,X)/\delta^{(\lambda)}(c),\dots,
u^{(\lambda,K)}(c,X)/\delta^{(\lambda)}(c)\in R[X].$$
\end{cor}

We leave it to the reader to formulate a similar result for B\'ezout domains,
using Corollary~\ref{DefBezout-Cor}.

\subsection*{Extremely coherent rings}
In the following lemma,
let $\cal C$ be a class of rings which is closed under direct products
and  $(\alpha,\beta)$-super coherent, and let $\gamma\colon\N^3\to\N$ be as
in Lemma~\ref{SupercoherentIsStablycoherent}.

\begin{lemma}
There exists a function
$\delta\colon\N^4\to\N$
with the following property: Given $R\in\cal C$ and $f_1,\dots,f_n\in R[X]$
of degree $\leq d$,
where $X=(X_1,\dots,X_N)$, there exist
solutions $y^{(1)},\dots,y^{(\gamma)}\in R[X]^n$
of degree $\leq\beta(N,d,n)$ to the homogeneous linear equation
\begin{equation}\label{HomEq}
f_1y_1+\cdots+f_ny_n=0
\end{equation}
in $R[X]^n$ such that every solution $y\in R[X]^n$ to \eqref{HomEq}
can be written as
$$y=a_1y^{(1)}+\cdots+a_\gamma y^{(\gamma)}\qquad
(\gamma=\gamma(N,d,n))$$
with $a_1,\dots,a_\gamma\in R[X]$ of degree $\leq\delta(N,d,e,n)$,
where $e=\deg(y)$.
\end{lemma}
\begin{proof}
In order to establish the existence of a function $\delta$
with the required property, 
we first note that it suffices to show
the following seemingly weaker statement:
\begin{itemize}
\item[($\ast$)] For any $N,d,e,n\in\N$, $n\geq 1$, there exists an integer
$\delta=\delta(N,d,e,n)\geq 0$ such that
given $R\in\cal C$ and $f_1,\dots,f_n\in R[X]=R[X_1,\dots,X_N]$
of degree $\leq d$, where $X=(X_1,\dots,X_N)$, 
every solution $y\in R[X]^n$ to \eqref{HomEq}
of degree $\leq e$ can be written as  
$$
y=b_1z^{(1)}+\cdots+b_\gamma z^{(\gamma)}
$$
with certain solutions $z^{(1)},\dots,z^{(\gamma)}\in R[X]^n$ of degree
$\leq\beta(N,d,n)$ and $b_1,\dots,b_\gamma\in R[X]$ of degree $\leq\delta$.
\end{itemize}
For suppose we have established this statement, and let 
$R\in\cal C$ and $f_1,\dots,f_n\in R[X]=R[X_1,\dots,X_N]$ of degree $\leq d$
be given. By the proof of
Lemma~\ref{SupercoherentIsStablycoherent} 
there exist $\gamma=\gamma(N,d,n)$
generators $y^{(1)},\dots,y^{(\gamma)} \in R[X]^n$ 
for the $R$-module of solutions
to \eqref{HomEq} of degree $\leq\beta(N,d,n)$.
By ($\ast$) any solution $y\in R[X]^n$ to \eqref{HomEq} of degree $\leq e$ can
be written in the form $y=b_1z^{(1)}+\cdots+b_\gamma z^{(\gamma)}$ 
with certain solutions 
$z^{(1)},\dots,z^{(\gamma)}\in R[X]^n$ of degree
$\leq\beta(N,d,n)$ and $b_1,\dots,b_\gamma\in R[X]$ of degree $\leq\delta$.
Expressing each $z^{(i)}$ as an $R$-linear combination of the $y^{(j)}$ yields
$y$ as an $R[X]$-linear combination of $y^{(1)},\dots,y^{(\gamma)}$ 
with coefficients of degree $\leq\delta$, as required.

Now suppose for a contradiction that ($\ast$) is false, 
that is, there exist $N,d,e,n\in\N$, $n\geq 1$, such that for
every $k\in\N$ we find a ring $R^{(k)}\in\cal C$, $f_1^{(k)},\dots,f_n^{(k)}\in
R^{(k)}[X]=R^{(k)}[X_1,\dots,X_N]$ of degree $\leq d$, and a solution
$y^{(k)}\in R^{(k)}[X]^n$ to the equation
$f_1^{(k)}y_1+\cdots+f_n^{(k)}y_n=0$ of degree $\leq e$ which cannot
be written as a linear combination of $\gamma(N,d,n)$
solutions of degree $\leq\beta(N,d,n)$ with coefficients in $R^{(k)}[X]$ of
degree $\leq k$. Put $R^*=\prod_k R^{(k)}$ and $f_i=(f_i^{(k)})\in
R^*[X]$, $y_i=(y_i^{(k)})\in R^*[X]$. Then $y=[y_1,\dots,y_n]^{\tr}$ is a solution
to the homogeneous equation $f_1y_1+\cdots+f_ny_n=0$ which cannot be
written as a linear combination of $\gamma(N,d,n)$ solutions of degree
$\leq\beta(N,d,n)$. By virtue of Lemma~\ref{SupercoherentIsStablycoherent}, 
this contradicts the fact that $R^*\in\cal C$ is
$(\alpha,\beta)$-coherent.
\end{proof}

\begin{remark}
The proof shows that if there exists an integer $k>0$ such that every
$R\in\cal C$ is of finite rank $\leq k$, then the function
$(N,d,e,n)\mapsto\delta(N,d,e,n)$  
can be chosen not to depend on $n$. 
\end{remark}

Let us call a class ${\cal C}$ of rings 
{\bf $(\alpha,\beta,\delta)$-extremely coherent} if $\cal C$ is
$(\alpha,\beta)$-super coherent and 
$\delta\colon\N^4\to\N$ satisfies the conclusion of the
previous lemma, with $\gamma$ as in 
Lemma~\ref{SupercoherentIsStablycoherent}. 
We say that $\cal C$ is extremely coherent if it is
$(\alpha,\beta,\delta)$-extremely coherent for some choice of
{\bf uniformity functions} $\alpha,\beta,\delta$.
The last lemma yields 
a refinement of the ``if'' direction of Theorem~\ref{SuperIsStable}:

\begin{cor}\label{SuperIsStable-Cor}
A class of rings which is closed under direct products and stably coherent
is extremely coherent.
\end{cor}

We say that 
a ring $R$ is $(\alpha,\beta,\delta)$-extremely coherent if the class
${\cal C}=\{R\}$ is $(\alpha,\beta,\delta)$-extremely coherent.

\begin{lemma}
Let $\alpha$, $\beta$, $\gamma$ and $\delta$ be as above.
\begin{enumerate}
\item Let $\{R^{(k)}\}_{k\in\N}$ be a family of 
$(\alpha,\beta,\delta)$-extremely coherent rings. 
For any filter $\cal F$ on $\N$, the reduced product
$\prod_k R^{(k)}/\cal F$ is $(\alpha,\beta,\delta)$-extremely 
coherent.
\item Let $R\subseteq S$ be a faithfully flat ring extension. If
$R$ is $\alpha$-uniformly coherent and
$S$ is $(\alpha,\beta,\delta)$-extremely coherent, then $R$
is $(\alpha,\beta,\delta)$-extremely coherent. 
\end{enumerate}
\end{lemma}
\begin{proof}
Part (1) follows from Theorem~\ref{Horn} on Horn formulas.
For (2) suppose that $S$ is an 
$(\alpha,\beta,\delta)$-extremely coherent ring,
faithfully flat over the $\alpha$-uniformly coherent subring $R$. 
Then $R$ is $(\alpha,\beta)$-super coherent, see the remarks
following Definition~\ref{Definition-Supercoherence}.
Let $f_1,\dots,f_n\in R[X]$ be of degree $\leq d$.
By Lemma~\ref{SupercoherentIsStablycoherent} there exist
generators  $y^{(1)},\dots,y^{(\gamma)}\in R[X]^n$ of 
degree $\leq\beta(N,d,n)$ for the module of syzygies of
$f=(f_1,\dots,f_n)$ in $R[X]$; here  $\gamma=\gamma(N,d,n)$.
Since $S$ is $(\alpha,\beta,\delta)$-extremely coherent, there also exist
syzygies $z^{(1)},\dots,z^{(\gamma)}\in S[X]^n$
of $f$ of degree $\leq\beta(N,d,n)$ 
such that every syzygy $y\in S[X]^n$ of $f$
can be written as a linear combination
$$y=b_1z^{(1)}+\cdots+b_\gamma z^{(\gamma)}$$
with $b_1,\dots,b_\gamma\in S[X]$ of degree 
$\leq\delta(N,d,e,n)$, where $\gamma=\gamma(N,d,n)$ and $e=\deg(y)$.
By faithful
flatness of $S$ over $R$, every $z^{(j)}$ is an $S$-linear combination
of $y^{(1)},\dots,y^{(\gamma)}$, and if we have $y\in R[X]^n$, 
then $y=a_1y^{(1)}+\cdots+a_{\gamma}y^{(\gamma)}$ for
some $a_1,\dots,a_{\gamma}\in R[X]$ of degree
$\leq \delta(N,d,e,n)$, where $\gamma=\gamma(N,d,n)$.
\end{proof}

\begin{question}
Is the class of $(\alpha,\beta)$-super coherent rings closed under
direct products? Equivalently, is there $\delta\colon\N^4\to\N$
such that every $(\alpha,\beta)$-super coherent
ring is $(\alpha,\beta,\delta)$-extremely coherent?
(The equivalence follows from part (1) of the previous lemma and
Corollary~\ref{SuperIsStable-Cor} applied to
${\cal C}=$ the class of all
$(\alpha,\beta)$-coherent rings.)
\end{question}

Before we give a list of examples and further questions, 
let us remark that if $R$ is a coherent
ring, then the canonical embedding of $R$ into the direct product
$S=\prod_{\frak m} R_{\frak m}$ of its localizations $R_{\frak m}$ at
maximal ideals $\frak m$ of $R$ makes $S$ into a faithfully flat $R$-module.
(See, e.g., \cite{Sabbagh}, p.~502--503.) 

\begin{examples}\ 

\begin{enumerate}
\item By the results of Hermann \cite{Hermann} and Seidenberg
\cite{Seidenberg1}, the class of fields is extremely coherent with uniformity 
functions $\alpha(n)=n$,
$\beta(N,d)=(2d)^{2^{N-1}}$  and $\delta(N,d,e)=(2d')^{2^{N-1}}$ for $N>0$, 
where
$d'=\max\bigl\{e,\beta(N,d)\bigr\}$.
\item The localization $R_{\frak m}$ of a von Neumann regular ring 
$R$ at one of its maximal ideals $\frak m$ is a field.
By the remark above, the previous lemma and example (1), this 
implies that the class of von Neumann regular
rings is extremely coherent with the same uniformity functions as in (1).
(This was first observed by Sabbagh \cite{Sabbagh}.)
\item The class of DVRs is extremely coherent,
by the proof of Theorem~4.1 in \cite{maschenb-ideal2}; see also
the remark following Corollary~\ref{SemiIsSuper}.
\item The class of hereditary rings is extremely coherent
with the same uniformity functions as in example (3),
by the remark above and the last lemma.
\item The class $\cal S$ of semihereditary rings is extremely coherent, by
Corollary~\ref{SuperIsStable-Cor}. The nature of the associated
uniformity functions $\beta$ and $\delta$ is somewhat mysterious. Can they
be chosen to be {\it doubly exponential}\/ similar to the ones in example (1)?
(In trying to answer this question it is enough
to restrict to the 
subclass of $\cal S$ consisting of all {\it valuation rings.}\/)
\end{enumerate}
\end{examples}

\subsection*{Rings with nilpotents}
So far, we have concentrated on classes of {\it reduced}\/ rings,
such as the class $\cal S$ of semihereditary rings. We will now
exhibit certain extremely coherent classes of rings extending $\cal S$
which also contain rings with non-zero nilradical.
We say that a module $M$ over a ring $R$
is {\bf $m$-presented} (for a given $m\geq 1$) if there exists an
exact sequence $R^m \to R^m \to M\to 0$ of $R$-linear maps. 
Using Theorem~\ref{Horn} it is routine to show:
%The minimal number of generators of a finitely
%generated $R$-module $M$ is given by $\dim_\kappa M\otimes_R \kappa$,
%hence bounded by the length of $M$.

\begin{lemma}\label{m-presented}
Let $\{R^{(k)}\}_{k\in\N}$ be a family of rings and
for each $k$ let $M^{(k)}$ be an $m$-presented $R^{(k)}$-module. Then
$M=\prod_k M^{(k)}$ is an $m$-presented $R$-module, where $R=\prod_k R^{(k)}$.
\end{lemma}

\begin{example}
If $R$ is a local Noetherian ring, 
then any finitely generated $R$-module of length $m$ is $m$-presented.
(An easy consequence of Nakayama's Lemma.)
\end{example}

%\begin{proof}
%The direct product $\prod_k (R^{(k)})^m$ is a module over $R=\prod_k R^{(k)}$,
%and we have a natural isomorphism $R^m = \left(\prod_k R^{(k)}\right)^m \iso
%\prod_k (R^{(k)})^m$ of $R$-modules. Hence the exact sequences 
%$$(R^{(k)})^m \to (R^{(k)})^m \to M^{(k)}\to 0 \qquad (k\in\N)$$ 
%give rise to an exact
%sequence $$R^m\to R^m \to M=\prod_k M^{(k)}\to 0$$ as required.
%\end{proof}

\begin{definition} 
For fixed $m\geq 1$ let ${\cal S}_m$ be the class of rings $R$ such that
\begin{enumerate}
\item $\Nil(R)$ is nilpotent of index $\leq m$;
\item $\Nil(R)$ is $m$-presented;
\item $R/\Nil(R)$ is semihereditary.
\end{enumerate}
Clearly we have ${\cal S}={\cal S}_1\subseteq {\cal S}_m$ for every $m$.
\end{definition}

\begin{prop}\label{Cm-Prop}
The class ${\cal S}_m$ is closed under direct products and extremely coherent.
It contains all Artinian rings of length $\leq m$.
\end{prop}
\begin{proof}
Let $R^{(k)}$ ($k\in\N$) be rings whose nilradical is nilpotent of 
index $\leq m$. Then the same is true for $R=\prod_k R^{(k)}$,
and $\Nil(R)=\prod_k\Nil(R^{(k)})$. Moreover,
if each $\Nil(R^{(k)})$ is an $m$-presented $R^{(k)}$-module, 
then $\Nil(R)$ is an $m$-presented $R$-module, by Lemma~\ref{m-presented}.
If each quotient ring $R^{(k)}/\Nil(R^{(k)})$ is semihereditary, then so is
$R/\Nil(R)\iso \prod_k R^{(k)}/\Nil(R^{(k)})$, by Lemma~\ref{Product}.
It follows that ${\cal S}_m$ is closed under direct products.

In order to show that ${\cal S}_m$ is extremely coherent, it remains to show
(by Corollary~\ref{SuperIsStable-Cor})
that for every $R\in{\cal S}_m$ and
every integer $N\geq 0$ the polynomial ring $R[X]=R[X_1,\dots,X_N]$ is
coherent. Since $R/\Nil(R)$ is semihereditary, the ring 
$(R/\Nil(R))[X]$ is coherent,
by Theorem~\ref{Vasconcelos}.
The natural surjection $R\to R/\Nil(R)$ induces a ring homomorphism
$$\phi\colon R[X]\to (R/\Nil(R))[X]$$ with finitely generated nilpotent kernel
$\ker\phi=\Nil(R)R[X]$. Now $\Nil(R)$ is a finitely presented $R$-module, 
hence $\Nil(R)R[X]$ is a finitely presented $R[X]$-module. (Since
the ring extension $R\to R[X]$ is faithfully flat.) It follows that  
$R[X]$ is coherent, by part (2) of Proposition~\ref{Descent}.

Every Artinian ring is isomorphic to a finite direct product of
local Artinian rings. Hence, since ${\cal S}_m$ is closed under direct
products, it suffices to show that ${\cal S}_m$ contains all local
Artinian rings $(R,{\frak m})$ of length $\leq m$. In this case we have
$\Nil(R)=\frak m$ and ${\frak m}^m=\{0\}$. By the remark before 
Lemma~\ref{m-presented}, $\frak m$ is $m$-presented. Moreover
$R/\Nil(R)=R/\frak m$ is a field, hence semihereditary. 
Therefore $R\in{\cal S}_m$.
\end{proof}

\begin{cor}\label{Artinian-Cor}
For each triple $(N,d,l)\in\N^3$ there exists a natural
number $\beta=\beta(N,d,l)$ such that for every Artinian ring $R$ of length
at most $l$ and polynomials $f_1,\dots,f_n\in R[X]=R[X_1,\dots,X_N]$ of
degree at most $d$, the module of syzygies of $(f_1,\dots,f_n)$ in $R[X]$ can
be generated by elements of $R[X]^n$ of degree at most $\beta$. 
\end{cor}

\begin{remark}
For Artinian local rings, the previous corollary was first 
proved by Schou\-tens \cite{Schoutens}.
\end{remark}

\section{Inhomogeneous Linear Equations in Polynomial Rings}\label{Inhomogeneous-Section}

The main purpose of this section is to show Theorems B and C from the
Introduction. On our way to proving Theorem~C we will also treat the
question of defining membership in the nilradical of a
finitely generated ideal in a polynomial ring over an arbitrary ring.

\subsection*{Uniform rings}
Let $\beta\colon\N^3\to\N$.
We say that a ring $R$ is {\bf $\beta$-uniform} if 
for given $N,d,n\in\N$ and
$f_1,\dots,f_n\in R[X]=R[X_1,\dots,X_N]$ of degree at most $d$, 
if $1\in (f_1,\dots,f_n)R[X]$, then there exist
$g_1,\dots,g_n\in R[X]$ of degree at most $\beta(N,d,n)$ such that
$$1=f_1g_1+\cdots+f_ng_n.$$
We say that $R$ is {\bf uniform} if $R$ is $\beta$-uniform
for some function $\beta$ as above.
A class $\cal C$ of rings is called $\beta$-uniform 
if every $R\in\cal C$ is $\beta$-uniform, and we say that
$\cal C$ is uniform if $\cal C$ is $\beta$-uniform for some $\beta$.

If a ring $R$ is uniform and of finite rank, then $R$ is 
$\beta$-uniform for some function $(N,d,n)\mapsto
\beta(N,d,n)$ which does not depend on $n$. (See the remark following
Lemma~\ref{SupercoherentIsStablycoherent}.)
The proof of the next lemma is similar to the proof of Lemma~\ref{SuperCoh}. 
We leave the details to the reader.

\begin{lemma}\label{Uniform-Lemma}
Let $\cal C$ be a class of rings.
The following are equivalent:
\begin{enumerate}
\item $\cal C$ is uniform.
\item For every family $\{R^{(k)}\}_{k\in\N}$ of rings in $\cal C$, every
filter $\cal F$ on $\N$ and
every $N\in\N$, if $I$ is a finitely generated ideal of 
$R^*[X]$ with $1\notin I$, then $1\notin I R[X]^*$, where $X=(X_1,\dots,X_N)$.
\item For every family $\{R^{(k)}\}_{k\in\N}$ of rings in $\cal C$ and every
$N\in\N$, a finitely generated proper ideal of
$\left(\prod_k R^{(k)}\right)[X]$ remains proper after extension to
$\prod_k R^{(k)}[X]$, with $X=(X_1,\dots,X_N)$.
\end{enumerate}
\end{lemma}

Of particular interest are classes of rings which are both uniform
and super coherent. Standard arguments (see, e.g., 
\cite{maschenb-ideal2}) show:

\begin{lemma}\label{Matrix-Lemma}
Let $\cal C$ be a uniform and super coherent class of rings. Then,
for any given $N,d,m,n\in\N$, $m,n\geq 1$,
there exists a natural number $\beta_m=\beta_m(N,d,n)$ with the
property that for all $R\in\cal C$, all
$m\times n$-matrices $A$ and all column vectors $b$
with entries in $R[X]=R[X_1,\dots,X_N]$ of degree at most $d$:
if the system $$Ay=b$$ is solvable in $R[X]$, then it has a solution in $R[X]$
all of whose entries have degree bounded from above by $\beta_m$. 
\rom{(In particular, $\cal C$ is extremely coherent.)}
\end{lemma}

Using the fact that a ring extension $R\subseteq S$ is faithfully flat
if and only if $S$ is a flat $R$-module and
$IS\neq S$ for every finitely generated ideal $I\neq R$ of
$R$, Lemmas \ref{SuperCoh} and \ref{Uniform-Lemma} imply:

\begin{cor}
Let $\cal C$ be a uniformly coherent class of rings. The following are
equivalent:
\begin{enumerate}
\item $\cal C$ is uniform and super coherent.
\item For every family $\{R^{(k)}\}_{k\in\N}$ of rings in $\cal C$, every
filter $\cal F$ on $\N$ and
every $N\in\N$, $R[X]^*$ is faithfully flat over
$R^*[X]$, where $X=(X_1,\dots,X_N)$.
\item For every family $\{R^{(k)}\}_{k\in\N}$ of rings in $\cal C$ and every
integer $N\geq 0$, the rings
$\prod_k R^{(k)}[X]$ is faithfully 
flat over its subring
$\left(\prod_k R^{(k)}\right)[X]$, with $X=(X_1,\dots,X_N)$.
\end{enumerate}
\end{cor}

The first goal of this section is to show that
uniformity is a rather serious restriction:

\begin{theorem}\label{Uniform}
A class $\cal C$  of rings is uniform if and only if there exists $m\geq 1$
such that $\Nil(R)$ is nilpotent of nilpotency index $\leq m$ and $R/\Nil(R)$
is von Neumann regular, for every $R\in\cal C$.
\end{theorem}

Before we begin the proof, we establish several lemmas.
For a proof of the first one see \cite{Cherlin}, Lemma~2.3.
An element $r$ of a ring $R$ is called {\bf regular} if
$r^2$ divides $r$ in $R$. If $r$ is both regular and nilpotent, then $r=0$.

\begin{lemma}\label{Cherlin-Lemma}
Let $R$ be a ring, $r\in R$, $m\in\N$. Then $r^{m+1}$ divides $r^m$
if and only if $r=r_1+s$ with $r_1\in R$ regular and $s\in R$, $s^m=0$.
\end{lemma}

The next lemma shows that the class of von Neumann regular
rings is uniform:

\begin{lemma}\label{vNregular}
Let $R$ be a von Neumann regular ring, $N\geq 0$, and $f_1,\dots,f_n\in
R[X]=R[X_1,\dots,X_N]$ of degree $\leq d$. If $1\in (f_1,\dots,f_n)R[X]$, then
there exist polynomials 
$g_1,\dots,g_n\in R[X]$ of degree $\leq \beta(N,d)=d^{N+1}$
such that
$$1=f_1g_1+\cdots+f_ng_n.$$
\end{lemma}
\begin{proof}
Since $R$ is von Neumann regular, $R$ can be embedded into a direct product
$S=\prod_{i\in I} K_i$ of a family of fields with $S$ faithfully
flat over $R$. Hence $S[X]$ is faithfully flat over $R[X]$. 
Replacing $R$ by $S$ if necessary we can therefore assume that 
$R$ is a direct product of a family of
fields, and in this case the lemma follows from the effective Nullstellensatz
of Koll\'ar \cite{Kollar}.
\end{proof}

Recall the familiar multinomial formula: For $e,M\in\N$, $M\geq 1$
\begin{equation}\label{Multinomial}
(Y_1+\cdots+Y_M)^e = \sum_{e_1+\cdots+e_M=e} \binom{e}{e_1,\dots, e_M}
Y_1^{e_1}\cdots Y_M^{e_M},
\end{equation}
where $Y_1,\dots,Y_M$ are distinct indeterminates over $\Z$ and
$\binom{e}{e_1,\dots, e_M}= \frac{e!}{e_1!\cdots e_M!}$ for all
$(e_1,\dots,e_M)\in\N^M$ with $e_1+\cdots+e_M=e$.
We record the following immediate consequence:

\begin{lemma}\label{Multinomial-Lemma}
Let $R$ be a ring whose nilradical $\Nil(R)$ is nilpotent of index $m$.
Then the nilradical $\Nil(R[X])$ of $R[X]$ 
is nilpotent of index $\leq \binom{N+d}{N}\cdot m$.
\end{lemma}

We can now prove Theorem~\ref{Uniform}.
The ``only if'' direction is implicit in the proof of Proposition~5 in
\cite{Sabbagh}. For the convenience of the reader 
we repeat the argument: Suppose that $\cal C$ is $\beta$-uniform,
let $R\in\cal C$ and $r\in R$ be arbitrary, $n\in\N$, 
and consider the following elements of $R[X]$
(where $X$ is a single indeterminate): $P(X)=rX+1$, $P_n(X)=r^n$.
Then obviously $P_n = r^nP-XP_{n+1}$, hence $1\in (P,P_n)$. Put
$m:=\beta(1,1,2)+1$, so there are polynomials $Q(X),Q_{m}(X)\in R[X]$
of degree $< m$ such that $1=PQ+P_{m}Q_{m}$. A computation now shows
that $r^{m+1}$ divides $r^{m}$. Hence by Lemma~\ref{Cherlin-Lemma}
$r=r_1+s$ with $r_1$ regular and $s^{m}=0$. Hence $\Nil(R)$ is nilpotent
of index $\leq m$, and $R/\Nil(R)$ is von Neumann regular. 

Conversely, suppose that there exists some $m\geq 1$ such that for
every $R\in\cal C$, $\Nil(R)$ is nilpotent of index $\leq m$ and $R/\Nil(R)$ is
von Neumann regular. Let $R\in\cal C$, $N\geq 1$ and
$f_1,\dots,f_n\in R[X_1,\dots,X_N]=R[X]$ of degree $\leq d$ with
$1\in (f_1,\dots,f_n)R[X]$. 
Hence $1\in (\bar{f_1},\dots,\bar{f_n})R/\Nil(R)[X]$,
where $\bar{f}$ denotes the image of the polynomial $f\in R[X]$ 
under the canonical surjection $R[X]\to (R/\Nil(R))[X]$.
By Lemma~\ref{vNregular}, there exist
$g_1,\dots,g_n\in R[X]$ of degree $\leq d^{N+1}$ such that
$$h := 1 + f_1g_1+\cdots+f_ng_n \in \Nil(R)R[X].$$
The degree of $h$ is at most $d^{N+2}$.
By Lemma~\ref{Multinomial-Lemma} 
it follows that $h^D=0$, where $D=\binom{N+d^{N+2}}{N}\cdot 
m$. On the other hand 
we have, by letting $M=n+1$ in \eqref{Multinomial} and
specializing $Y_1$ to $1$ and $Y_2,\dots,Y_M$ to $f_1g_1,\dots,f_ng_n$,
respectively:
$$h^D = 1 - (f_1h_1+\cdots+f_nh_n)$$
with $h_1,\dots,h_n\in R[X]$ of degree $\leq Dd^{N+2}$.
Hence $R$ is $\beta$-uniform with $$\beta(N,d)=\binom{N+d^{N+2}}{N}\cdot md^{N+2}.$$
This finishes the proof of Theorem~\ref{Uniform}.
(Note that $\beta$ does not depend on $n$ and
is even {\it linear}\/ in the upper bound $m$ on the 
nilpotency index.) \qed

\begin{cor}
Let $R$ be a ring.
\begin{enumerate}
\item If $\Nil(R)$ is finitely generated, then $R$ is
uniform and super coherent if and only if $\Nil(R)$ is finitely presented
and $R/\Nil(R)$ is von Neumann regular. 
\item If $R$ is Noetherian, then $R$ is uniform 
if and only if 
$R/\Nil(R)$ is semisimple \rom{(}i.e., isomorphic to 
a finite direct product of fields\rom{)}. In particular, 
a Noetherian uniform ring is super
coherent.
\end{enumerate}
\end{cor}
\begin{proof}
For the first part use Proposition~\ref{Cm-Prop}; for the second part note
that the semisimple rings are exactly the Noetherian von
Neumann regular rings.
\end{proof}

Combining Theorem~\ref{Uniform} with Corollary~\ref{Artinian-Cor} yields
the following result:

\begin{cor}\label{Artinian-Cor-2}
For each triple $(N,d,l)\in\N^3$ there exists 
$\beta=\beta(N,d,l)\in\N$ such that for every Artinian ring $R$ of length
at most $l$ and polynomials $f_0,f_1,\dots,f_n\in R[X]=R[X_1,\dots,X_N]$ of
degree at most $d$: if $f_0\in (f_1,\dots,f_n)R[X]$, then 
$$f_0=f_1g_1+\cdots+f_ng_n$$
for some $g_1,\dots,g_n\in R[X]$ of degree at most $\beta$.
\end{cor}
\begin{proof}
Let $R$ be an Artinian ring of length $\leq l$, and
$f_0,f_1,\dots,f_n\in R[X]=R[X_1,\dots,X_N]$ of
degree at most $d$.
Consider the homogeneous linear equation
\begin{equation}\label{Artinian-Cor-Eq}
f_0y_0+f_1y_1+\cdots+f_ny_n=0.
\end{equation}
By Corollary~\ref{Artinian-Cor} we find
generators $y^{(1)},\dots,y^{(K)}\in R[X]^{n+1}$ for the
module of solutions to \eqref{Artinian-Cor-Eq} whose degrees are
uniformly bounded in terms of $N$, $d$, $l$ (independent of $R$ and
$f_0,\dots,f_n$). For $g_1,\dots,g_n\in R[X]$ 
we have $f_0=f_1g_1+\cdots+f_ng_n$ if and only if
$(1,-g_1,\dots,-g_n)^{\tr}$ is a solution to \eqref{Artinian-Cor-Eq}.
Write
$y^{(k)}=\bigl(y^{(k)}_0,\dots,y^{(k)}_n\bigr)^{\tr}$.
By Theorem~\ref{Uniform}, if $1\in\bigl(y^{(1)}_0,\dots,y^{(K)}_0\bigr)R[X]$,
then there exist $h_1,\dots,h_K\in R[X]$ with
$1=y^{(1)}_0h_1+\cdots+y^{(K)}_0h_K$ whose degrees are uniformly
bounded in terms  of $N$, $d$ and $l$.
The corollary follows.
\end{proof}

\begin{remark}
The last corollary was first 
proved by Schou\-tens \cite{Schoutens} (for local Artinian rings).
For $l=1$ and $R$ local, we recover Hermann's theorem quoted after Theorem~A.
\end{remark}

We now turn to issues of definability. 

\subsection*{Definability of membership in the nilradical}

In the rest of this section we 
let $C=(C_1,\dots,C_M)$ be a tuple of parametric variables. Let
$$f_0(C,X),f_1(C,X),\dots,f_n(C,X)\in\Z[C,X].$$ 
%Recall that
%given a ring $R$, we write $\sqrt{I}$ 
%for the nilradical of an ideal $I$ of $R$. (See \S\ref{NotationsSection}.)
For any field $K$ and $c\in K^M$, we have
\begin{equation}\label{Nil1}
f_0(c,X)\in\sqrt{\bigl(f_1(c,X),\dots,f_n(c,X)\bigr)K[X]}
\quad \Longleftrightarrow
\end{equation}
$$\text{for all $a\in (K^{\alg})^N$: 
$\bigl(f_1(c,a)=0\wedge\cdots\wedge f_n(c,a)=0\bigr) 
\Rightarrow f_0(c,a)=0$},$$
by Hilbert's Nullstellensatz. (Here $K^{\alg}$ denotes an 
algebraic closure of $K$.)
Hence, using primitive recursive
quantifier elimination for the theory of algebraically closed fields,
we may find, primitive recursively in $f_0,\dots,f_n$, a family
$(p_{ij},q_i)_{\substack{1\leq i\leq m \\ 1\leq j\leq k}}$ with $k\in\N$,
$p_{ij}(C)\in\Z[C]$, $q_i\in\Z[C]$, such that for all fields $K$ and
$c\in K^M$,
\begin{equation}\label{Nil2}
f_0(c,X)\in \sqrt{\bigl(f_1(c,X),\dots,f_n(c,X)\bigr)K[X]}
\quad\Longleftrightarrow
\end{equation}
$$\bigwedge_{i=1}^m \bigl(p_{i1}(c)=0\wedge\cdots\wedge p_{ik}(c)=0
\Rightarrow q_i(c)=0\bigr).$$
In other words, we have for all fields $K$ and $c\in K^M$,
$$f_0(c,X)\in \sqrt{\bigl(f_1(c,X),\dots,f_n(c,X)\bigr)K[X]}
\quad\Longleftrightarrow$$
$$q_i(c)\in\sqrt{\bigl(p_{i1}(c),\dots,p_{ik}(c)\bigr)K}
\quad\text{for all $i=1,\dots,m$.}$$
(Since in a field $K$, the nilradical of an ideal is either equal
to $K$ or to $(0)$.)
We now want to show that this equivalence in fact holds for {\it
all}\/ rings $R$ in place of $K$ and parameter tuples $c\in R^M$.
(This was pointed out to us by van den Dries.)

In the following let $R$ be an arbitrary ring. 
For $c\in R$ and ${\frak p}\in\Spec R$ we write $c/{\frak p}:=c+{\frak p}\in
R/{\frak p}$; more generally, if $c=(c_1,\dots,c_M)\in R^M$, then
we write $c/{\frak p}$ for 
$(c_1/{\frak p},\dots,c_M/{\frak p})\in (R/{\frak p})^M$.
For a polynomial $f\in
R[X]=R[X_1,\dots,X_N]$ and an ideal $I$ of $R[X]$ we denote by
$f_{(\frak p)}$ and $I_{(\frak p)}$ the image of $f$ and $I$,
respectively, under the canonical homomorphism
$$R[X]\to (R/\frak p)[X] \einb k_{\frak p}[X],$$
where ${\frak p}\in\Spec R$, $k_{\frak p}:=\Frac(R/{\frak p})$.

\begin{lemma}\label{NilLemma}
For $f\in R[X]$ and an ideal $I$ of $R[X]$, we have:
$$f\in\sqrt{I} \qquad\Longleftrightarrow\qquad \text{$f_{(\frak p)}\in
\sqrt{I_{(\frak p)}}\quad$ for all ${\frak p}\in\Spec R$ with 
${\frak p}\supseteq I\cap R$.}$$
\end{lemma}
\begin{proof}
The direction ``$\Rightarrow$'' is trivial. Suppose
$f\notin\sqrt{I}$. Then there exists a prime ideal ${\frak P}\supseteq I$
such that $f\notin\frak P$. Let ${\frak p}={\frak P}\cap R$, and let
$\bar{X_1},\dots,\bar{X_N}$ be the images of $X_1,\dots,X_N$ under the
canonical homomorphism $R[X]\to R[X]/{\frak P}=S$. We may naturally
identify $R/\frak p$ with a subring of $S$ and thus $k_{\frak p}$
with a subfield of $\Frac(S)$. We define a $k_{\frak p}$-homomorphism
$k_{\frak p}[X]\to\Frac(S)$ by $X_i\mapsto\bar{X_i}$ for
$i=1,\dots,N$. The image of $I_{(\frak p)}$ under this homomorphism is
$(0)$, so $(\bar{X_1},\dots,\bar{X_N})\in S^N$ is a zero of $I_{(\frak
p)}$, whereas the image of $f_{(\frak p)}$ is $0\neq f/{\frak P}\in
S$, so $(\bar{X_1},\dots,\bar{X_N})$ is not a zero of $f_{(\frak
p)}$. Thus $f_{(\frak p)}\notin\sqrt{I_{(\frak p)}}$.
\end{proof}

We now obtain the desired result:

\begin{prop}\label{NilProp}
For all $c\in R^M$, we have
$$f_0(c,X)\in \sqrt{\bigl(f_1(c,X),\dots,f_n(c,X)\bigr)R[X]}
\quad\Longleftrightarrow$$
$$q_i(c)\in\sqrt{\bigl(p_{i1}(c),\dots,p_{ik}(c)\bigr)R}
\quad\text{for all $i=1,\dots,m$.}
$$
\end{prop}
\begin{proof}
Suppose $f_0(c,X)\in \sqrt{\bigl(f_1(c,X),\dots,f_n(c,X)\bigr)R[X]}$, and let
$i\in\{1,\dots,m\}$ and
${\frak p}\in\Spec R$ with ${\frak
p}\supseteq\bigl(p_{i1}(c),\dots,p_{ik}(c)\bigr)R$. Then in $R/{\frak p}$, we have 
$$p_{i1}(c/{\frak p})=\cdots=p_{ik}(c/{\frak p})=0,$$
hence $q_i(c)\in\frak p$, and thus
$q_i(c)\in\sqrt{\bigl(p_{i1}(c),\dots,p_{ik}(c)\bigr)R}$. 
%Conversely,
%suppose $q_i(c)\in\sqrt{\bigl(p_{i1}(c),\dots,p_{ik}(c)\bigr)}$ for all
%$i=1,\dots,m$. 
Suppose that $$f_0(c,X)\notin
\sqrt{\bigl(f_1(c,X),\dots,f_n(c,X)\bigr)R[X]}.$$ Then there exists ${\frak
p}\in\Spec R$ such that $$f_0(c,X)_{(\frak p)}\notin 
\sqrt{\bigl(f_1(c,X)_{(\frak p)},\dots,f_n(c,X)_{(\frak p)}\bigr)k_{\frak p}[X]},$$
by the lemma;
thus for some $i\in\{1,\dots,m\}$, we have
$$p_{i1}(c/{\frak p})=\cdots=p_{ik}(c/{\frak p})=0, \quad q_i(c/{\frak
p})\neq 0,$$
by \eqref{Nil1}.
Therefore $q_i(c)\notin \sqrt{\bigl(p_{i1}(c),\dots,p_{ik}(c)\bigr)R}$.
\end{proof}

\begin{remark}
The ideal $c(f)$ of $R$ generated by
the coefficients of a polynomial $f\in R[X]$ is called 
the {\bf content} of $f$.
Lemma~\ref{NilLemma} may also be used to obtain a quick proof of the
following generalization of Gauss' Lemma: $$\sqrt{c(fg)}=
\sqrt{c(f)}\cdot\sqrt{c(g)}\qquad\text{for all $f,g\in R[X]$.}$$ To see this,
note first that it suffices to show the inclusion $\supseteq$. Moreover, it
is enough treat the case where $R=\Z$ and
the coefficients of $f$ and $g$ are pairwise distinct 
indeterminates over $\Z[X]$. Fixing an enumeration $X^{\mu_1},\dots,X^{\mu_M}$ 
of all monomials of degree $\leq d$ we may therefore write 
$f(C,X)=\sum_i C_iX^{\mu_i}\in\Z[C,X]$ and 
$g(C',X)=\sum_i C_i'X^{\mu_i}\in\Z[C',X]$ where
$C=(C_1,\dots,C_M)$ and $C'=(C_1',\dots,C_M')$; here $d\in\N$ and
$M=\binom{N+d}{N}$.
By Lemma~\ref{NilLemma} (applied to $\Z[C,C']$ in place of $R[X]$) 
we may further reduce to the case where $R=K$ is a prime field
(i.e., $K=\Q$ or $K=\F_p$ for some prime $p$). 
Since for all $c=(c_1,\dots,c_M),c'=(c_1',\dots,c_M')\in 
(K^{\alg})^{M}$ we have
$$f(c,X)\cdot g(c',X)=0 \quad\Longleftrightarrow\quad f(c,X)=0\text{ or }
g(c',X)=0,$$ the algebraic subset $V$ of
$(K^{\alg})^{2M}$ defined by the vanishing of the
coefficients of $f(C,X)\cdot g(C',X)$
is the union 
$$V=\bigl\{(c,c') : c_1=\cdots=c_M=0\bigr\}\cup 
\bigl\{(c,c') : c_1'=\cdots=c_M'=0\bigr\}.$$
The Nullstellensatz now 
yields the claim. (See \cite{Northcott} for a different proof.)
\end{remark}

Let ${\cal L}^*_{\rad}$ be the language of rings augmented by a $(k+1)$-ary
predicate symbol $\rad_k$, for each $k>0$. We construe a ring $R$
as an ${\cal L}^*_{\rad}$-structure by interpreting the ring symbols as usual
and the symbols $\rad_k$, for $k>0$, by
$$R\models\rad_k(r_0,r_1,\dots,r_k) \quad:\Longleftrightarrow\quad
r_0\in\sqrt{(r_1,\dots,r_k)R},$$
for $r_0,\dots,r_k\in R$. 

\begin{remark}
If $R$ is a B\'ezout domain and $r_0,r_1,\dots,r_k\in R$, then 
$$r_0\in\sqrt{(r_1,\dots,r_k)R}\quad  \Longleftrightarrow\quad 
r_0\in\sqrt{\gcd(r_1,\dots,r_k)R}.$$
In particular, if $R$ is a DVR with associated valuation $v$, then
$$r_0\in\sqrt{(r_1,\dots,r_k)R} \qquad\Longleftrightarrow\qquad
v(r_0)>0\vee\bigvee_{i=1}^k v(r_i)=0,$$
so the relations $\rad_k$ are quantifier-free definable in
the ${\cal L}_{\div}$-structure $R$.
\end{remark}

By the discussion above, we obtain:

\begin{cor}\label{CorRadical}
From the polynomials 
$f_0(C,X),\dots,f_n(C,X)\in\Z[C,X]$
one can primitive recursively
construct a quan\-ti\-fier-free ${\cal L}^*_{\rad}$-formula 
$\varphi(C)$ such that for all rings $R$ and all $c\in R^M$,
$$R\models\varphi(c) \quad\Longleftrightarrow\quad
f_0(c,X)\in\sqrt{\bigl(f_1(c,X),\dots,f_n(c,X)\bigr)R[X]}.$$
\end{cor}

In particular, from the polynomials $f_1(C,X),\dots,f_n(C,X)\in\Z[C,X]$  one can
primitive recursively construct a
quan\-ti\-fier-free ${\cal L}^*_{\rad}$-formula $\varphi(C)$ such that for 
every ring $R$, the set 
$$\bigl\{c\in R^M : 1\in\bigl(f_1(c,X),\dots,f_n(c,X)\bigr)R[X]\bigl\}$$
is defined by $\varphi$. 

\begin{remark}
The relation $\rad_1$ is indispensable for defining
membership in the nilradical of an ideal in $R[X]$ in a quantifier-free
way, as in the previous corollary. This can be shown by a
modification of the example in Section~6 of \cite{maschenb-ideal2}:
Let $a$, $b$ be elements of a ring $R$, and suppose that $X$ is a single
indeterminate. Then
$$1\in\bigl(1-aX,bX\bigr)R[X] \quad\Longleftrightarrow\quad
a\in \sqrt{bR}.$$
\end{remark}
\begin{proof}
If $a^n=bc$ for some $n\in\N$, $n>0$, and $c\in R$, then
$$1=(1+aX+\cdots+a^{n-1}X^{n-1})\cdot (1-aX)+cX^{n-1}\cdot bX,$$
exhibiting $1$ as an element of $\bigl(1-aX,bX\bigr)R[X]$.
Conversely, suppose that $$1\in \bigl(1-aX,bX\bigr)R[X].$$ Then $1-\bar{a}X$ 
is a unit in the ring $(R/bR)[X]$, where $\bar{a}=a+bR$. But in the
formal power series ring $(R/bR)[[X]]$, the element 
$1-\bar{a}X$ has multiplicative inverse
$$1+\bar{a}X+\bar{a}^2X^2+\bar{a}^3X^3+\cdots.$$ 
By uniqueness of inverses in $(R/bR)[[X]]$ it follows that 
$a\in \sqrt{bR}$ as required.
\end{proof}

Suppose $R$ is a computable ring such that for given elements $r_0,\dots,r_k$ of $R$ one can
decide whether $r_0\in\sqrt{(r_1,\dots,r_k)R}$. 
Then the computable ring $R[X]$ also has this
property, i.e., given $f_0,\dots,f_n\in R[X]$ one can 
effectively decide whether $f_0\in\sqrt{(f_1,\dots,f_n)R[X]}$.
For $R=\Z$, we have an even better result.
Namely, given $r_0,r_1,\dots,r_k\in\Z$, we can check in polynomial time
whether $r_0\in\sqrt{(r_1,\dots,r_k)\Z}$: we first find $a\in\Z$ such that
$(r_1,\dots,r_k)\Z=a\Z$, by the Euclidean Algorithm, and then we check
whether $a|(r_0)^e$, where $e$ is the integral part 
$\bigl[\log_2\abs{a}\bigr]$ if $a\neq 0$,
$e=1$ else. Thus validity of quantifier-free
${\cal L}^*_{\rad}$-formulas in $\Z$ can be checked in polynomial
time. This, together with the previous corollary, shows that for fixed
$f_0(C,X),\dots,f_n(C,X)\in\Z[C,X]$, membership in the set
$$\left\{c\in\Z^M : f_0(c,X)\in\sqrt{\bigl(f_1(c,X),\dots,f_n(c,X)\bigr)\Z[X]}\right\}$$
is decidable in polynomial time.
Moreover, we have a primitive
recursive algorithm which, upon input of $f_0,\dots,f_n\in\Z[X]$, 
decides whether $f_0\in\sqrt{(f_1,\dots,f_n)\Z[X]}$.

\subsection*{Definability of ideal membership}
We let %$C=(C_1,\dots,C_M)$ be a tuple of parametric variables,
$$A(C,X)=
\bigl(a_{ij}(C,X)\bigr)_{\substack{1\leq i\leq m\\
1\leq j\leq n}}$$ be an $m\times n$-matrix with entries $a_{ij}(C,X)
\in\Z[C,X]$, and $$b(C,X)=\begin{bmatrix} b_1(C,X) \\ \vdots \\ b_m(C,X)
\end{bmatrix}$$ with $b_i(C,X)\in\Z[C,X]$.

\begin{theorem}\label{SolvBezout}
There exists a quantifier-free ${\cal L}_{\rad}$-formula $\varphi(C)$
such that for all B\'ezout domains $R$, the set
\begin{equation}\label{SolvSet}
\bigl\{c\in R^M : \text{$A(c,X)y=b(c,X)$ is solvable in $R[X]$}\,\bigr\}
\end{equation}
is defined by $\varphi$.
\end{theorem}
\begin{proof}
Similar to the proof of Corollary~\ref{Artinian-Cor-2}, using
Corollary~\ref{DefBezout-Cor} and
the remarks following Corollary~\ref{CorRadical}.
\end{proof}

\begin{remark}
The case $m=1$, $R=\Z$ yields Theorem~C in the Introduction.
The remark about polynomial-time computability
after Theorem~C is a consequence of the discussion following
Corollary~\ref{CorRadical}.
\end{remark}

From the pair $(A,b)$
one can {\it construct}\/ (primitive recursively) a quantifier-free
${\cal L}_{\rad}$-formula $\varphi(C)$ which defines the set \eqref{SolvSet}
{\it in every PID $R$.}\/
Specializing even further to DVRs (and using
Corollary~\ref{DefBezout-Cor-2} 
instead of Corollary~\ref{DefBezout-Cor}) we get the following result.

\begin{cor}\label{SolvDVR}
From $(A,b)$ one can primitive recursively construct a 
quantifier-free ${\cal L}_{\div}$-formula $\psi(C)$
such that for all DVRs $R$, the set
$$
\bigl\{c\in R^M : \text{$A(c,X)y=b(c,X)$ is solvable in $R[X]$}\,\bigr\}
$$
is defined by $\psi$.
\end{cor} 

\begin{remarks}\ 

\begin{enumerate}
\item If $N=0$, then the quantifier-free formula $\varphi$ 
in Theorem~\ref{SolvBezout} may be chosen in the sublanguage
${\cal L}_{\gcd}$ of ${\cal L}_{\rad}$, and we can find a quantifier-free
formula $\psi$ which has the property in the previous corollary
{\it for all valuation rings $R$.}\/
This follows from \cite{vdDries-Holly} or more 
directly from a theorem of I.~Heger, 1856 (see \cite{OLeary-Vaaler}): 
if $R$ is a Pr\"ufer domain,
$A\in R^{m\times n}$ has rank $m$ and  $b\in R^m$ is a column vector, 
then $Ay=b$ has a solution $y\in R^n$ if and only
if the ideals of $R$ 
generated by all $m\times m$-minors of $A$ and 
by all $m\times m$-minors of $(A,b)$, respectively,
coincide.
\item Remark (1) remains true if more generally all polynomials 
$a_{ij}(C,X)$ are homogeneous in the indeterminates $X=(X_1,\dots,X_N)$. 
To see this, note that
for homogeneous $f_0,f_1,\dots,f_n\in R[X]$ 
with coefficients in a ring $R$ we have
$f_0\in (f_1,\dots,f_n)R[X]$ if and only if $f_0=f_1g_1+\cdots+f_ng_n$
for homogeneous polynomials $g_1,\dots,g_n\in R[X]$, with $g_j=0$ if 
$\deg f_j>\deg f_0$ and
$\deg g_j=\deg f_0-\deg f_j$ otherwise, for every $j$.
\end{enumerate}
\end{remarks}

\section{Prime Ideals}\label{Prime-Section}

Let $f_1(C,X),\dots,f_n(C,X)\in\Z[C,X]$, where again $C=(C_1,\dots,C_M)$.
In this last section, we want to apply the results obtained so far
to study instances of the following problem: Given a ring $R$, find
a description by a first-order formula (in a natural language)
of the set
\begin{equation}\label{PrimalitySet}
\bigl\{c\in R^M: \text{$\bigl(f_1(c,X),\dots,f_n(c,X)\bigr)R[X]$ 
is a prime ideal}\bigr\}.
\end{equation}
We first consider this question in the case that $R=K$ is a field.
By \cite{vdDries-Schmidt}, (2.10)~(ii), (iv), there exist
natural numbers $\alpha>1$ and $\beta$ (only depending on the $f_j$'s) 
such that for all fields $K$,
$c\in K^M$, and the ideal $I=\bigl(f_1(c,X),\dots,f_n(c,X)\bigr)$ of $K[X]$,
we have
$$\text{$I$ is radical} \qquad\Longleftrightarrow\qquad
\text{for all $f\in K[X]$ of degree $\leq\beta$:
$f^\alpha\in I \Rightarrow f\in I$}$$
and
\begin{multline*}
\text{$I$ is primary} \quad\Longleftrightarrow\quad\\
\text{$1\notin I$, and for all $f,g\in K[X]$ of degree $\leq\beta$: $fg\in I \Rightarrow f\in I$ or
$g^\alpha\in I$.}
\end{multline*}
In particular, there is a {\em universal}\/ formula in the language of rings
${\cal L}=\{{0},{1},{+},{-},{\cdot}\}$
defining the set of coefficients
$c\in K^M$ such that the ideal in $K[X]$ generated by $f_1(c,X),\dots,f_n(c,X)$ is radical, for
every field $K$; similarly for
``primary'' in place of ``radical''.
(If we restrict ourselves to algebraically closed $K$, then these
formulas may even be chosen quantifier-free, by quantifier-elimination of the theory of algebraically
closed fields.)
Since an ideal of a ring is prime if and only if it is radical and primary, we also get
\begin{multline*}
\text{$I$ is prime} \quad\Longleftrightarrow\quad\\
\text{$1\notin I$, and for all $f,g\in K[X]$ of degree $\leq\beta$: $fg\in I \Rightarrow f\in I$ or
$g\in I$,}
\end{multline*}
and there exists a universal ${\cal L}$-formula defining the set
\eqref{PrimalitySet} for all fields $R=K$. In \cite{vdDries-Thesis}, 
Chapter~IV, \S{}3, it was shown that
\eqref{PrimalitySet} may even be defined quantifier-free in a certain natural 
extension of ${\cal L}$, uniformly for all fields $R=K$. We give a brief
account of this result, simplifying it in the process by
replacing some of the Skolem functions
used in the extension of the language ${\cal L}$ by certain {\it 
predicate symbols}\/ for roots of separable
polynomials, and extending it to define the properties
{\it primary}\/ and {\it radical.}\/

\subsection*{Prime ideals in polynomial rings over fields}
Let $K$ be a field and $p=\operatorname{char} K$ if
$\operatorname{char} K>0$, $p=1$
if $\operatorname{char} K=0$. A field extension 
$L|K$ is called
\begin{enumerate}
\item {\bf separable} if $L^p$ and $K$ are linearly disjoint over
$K^p$,
\item {\bf primary} if the separable algebraic closure of $K$ in $L$
equals $K$, and
\item {\bf regular} if it is both separable and primary.
\end{enumerate}
The following lemma and its corollary
below are well-known:

\begin{lemma}
Let $A$ be a $K$-algebra and $B=A\otimes_K L$, an $L$-algebra.
\begin{enumerate}
\item If $L|K$ is separable and $A$ is reduced, then $B$ is reduced.
\item If $L|K$ is primary and $\Nil(A)$ is a prime ideal of $A$, then
$\Nil(B)$ is a prime ideal of $B$.
\item If $L|K$ is regular and $A$ is an integral domain, then $B$ is
an integral domain.
\end{enumerate}
\end{lemma}
\begin{proof}
Part~(1) follows from Proposition~5 of \cite{BA}, Chapitre~V, \S
15. For (2), note that replacing $A$ by $A/\Nil(A)$ we may assume that
$A$ is an integral domain. Now it follows from \cite{BA}, Chapitre~V,
\S 17, Corollaire~to Proposition~1, that $\Nil(B)$ is prime in
$B$. Since a ring is an integral domain if and only if it is reduced
and the set of its nilpotent elements is a prime ideal, (3)
follows from (1) and (2).
\end{proof}

\begin{cor}\label{ExtensionIdeal}
Let  $I$ be an ideal of $K[X]$.
\begin{enumerate}
\item If $L|K$ is separable and $I$ is a radical ideal, then $IL[X]$ is
a radical ideal.
\item If $L|K$ is primary and $I$ is a primary ideal, then $IL[X]$ is
a primary ideal.
\item If $L|K$ is regular and $I$ is a prime ideal, then $IL[X]$ is a
prime ideal. 
\end{enumerate}
\end{cor}

As remarked above, there exists a
quantifier-free ${\cal L}$-formula $\varphi(C)$ such that for all 
algebraically closed fields $K$ and $c\in K^M$:
\begin{equation}\label{RadicalEqu}
K\models\varphi(C) \quad\Longleftrightarrow\quad
\text{$\bigl(f_1(c,X),\dots,f_n(c,X)\bigr)K[X]$ is radical.}
\end{equation}
From (1) of the previous corollary, it follows immediately that the foregoing equivalence
\eqref{RadicalEqu} also holds for all {\em perfect}\/ fields $K$ and $c\in K^M$.

Let ${\cal L}_1$ be the language $\cal L$ 
of rings augmented by a unary function symbol $^{-1}$ and,
for every $m\geq 1$, an $m$-ary predicate symbol $Z_m$. We let $T_1$ be the extension of the theory of
rings (in the language $\cal L$ of rings) by the defining axiom
\begin{equation}\label{InvAxiom}
\forall x\bigl((x=0\wedge x^{-1}=0) \vee (x\neq 0 \wedge xx^{-1}=1)\bigr)
\end{equation}
and, for each $m\geq 1$, an axiom saying that for every model of $T_1$ with underlying field $K$ and
$(a_1,\dots,a_{m})\in K^m$,
\begin{multline*}
K\models Z_m(a_1,\dots,a_{m}) \quad\Longleftrightarrow\quad\\
\text{$T^m+a_{1}T^{m-1}+\cdots+a_m \in K[T]$ is separable and has a zero in $K$.}
\end{multline*}
Every field can be expanded uniquely to a model of $T_1$, and
a substructure of a model of $T_1$ is a field (but not necessarily a model of $T_1$).
Note that we include the symbol $^{-1}$ for convenience only: 
in $T_1$, every quantifier-free ${\cal L}_1$-formula is equivalent to a quantifier-free ${\cal L}_0$-formula,
where ${\cal L}_0={\cal L}_1\setminus\{{^{-1}}\}$.
%(One way to see this is to note that the
%theory $T_1$ has closures of ${\cal L}_0$-substructures
%in the sense of Chapter~\ref{PrelimChapter}, \S\ref{PrelimSection}.)

The following model-theoretic fact is proved by a standard application
of the Compactness Theorem; we leave the proof to the reader. 

\begin{lemma}\label{QE}
Let ${\cal L}$ and ${\cal L}^*$ be languages
\rom{(}in the sense of first-order logic\rom{)} with
${\cal L}\subseteq {\cal L}^*$, and let
$T^*$ be an ${\cal L}^*$-theory. For an ${\cal L}^*$-formula
$\varphi^*(x)$, $x=(x_1,\dots,x_n)$, the following are equivalent:
\begin{enumerate}
\item There exists a quantifier-free ${\cal L}$-formula $\varphi(x)$
such that $T^*\models\forall x(\varphi^*\leftrightarrow\varphi)$.
\item For all models ${\bf A}^*$ and ${\bf B}^*$ of $T^*$ whose reducts
to ${\cal L}$ have a common ${\cal L}$-substructure ${\bf C}=(C,\dots)$,
and for all $c\in C^n$: $${\bf A}^*\models\varphi^*(c) 
\quad\Longleftrightarrow\quad
{\bf B}^*\models\varphi^*(c).$$
\end{enumerate}
\end{lemma}

\begin{remark}
Suppose that one of the equivalent conditions in the lemma holds for an
${\cal L}^*$-formula $\varphi^*(x)$.
If ${\cal L}^*$ and $T^*$ are recursively enumerable, then a quantifier-free
${\cal L}$-formula $\varphi$ as in (1) can be found effectively, by
G\"odel's Completeness Theorem.
\end{remark}

For a field $K$, we denote the separable algebraic closure of $K$ (in a fixed algebraic closure of $K$)
by $K_{\sep}$. 

\begin{lemma}
Suppose $E$ and $F$ are the underlying fields of models of $T_1$ having a 
common ${\cal L}_1$-substructure with underlying field $K$. There exists an isomorphism
$$E\cap K_{\sep} \overset{\iso}{\longrightarrow} F\cap K_{\sep}$$ which is the identity on $K$.
\end{lemma}

This lemma is due to Ax (\cite{Ax}, \S{}3,
Lemma~5). We use it to show:

\begin{cor}\label{PrimaryCor}\label{PrimCor1}
There exists a quantifier-free ${\cal L}_0$-formula $\varphi_0(C)$ such that for every field $K$ and all
$c\in K^M$,
$$K\models\varphi_0(c)\quad\Longleftrightarrow\quad 
\text{$\bigl(f_1(c,X),\dots,f_n(c,X)\bigr)K[X]$ is primary.}$$
\end{cor}
\begin{proof}
By the discussion above, there exists an ${\cal L}$-formula $\varphi(C)$ (possibly involving quantifiers)
such that for all fields $K$ and $c\in K^M$,
$$K\models\varphi(c)\quad\Longleftrightarrow\quad 
\text{$\bigl(f_1(c,X),\dots,f_n(c,X)\bigr)K[X]$ is primary.}$$
Suppose now that $E$ and $F$ are the underlying fields of models of $T_1$ having a common
substructure with underlying field $K$, and suppose $c\in K^M$ is such that $E\models\varphi(c)$.
By the previous lemma,  it follows that
$F\cap K_{\sep} \models\varphi(c)$. Since the field extension $F\supseteq 
F\cap K_{\sep}$ is primary, we get
$F\models\varphi(c)$, by Corollary~\ref{ExtensionIdeal},~(2). 
By Lemma~\ref{QE},
$\varphi$ is equivalent to a quantifier-free ${\cal L}_1$-formula $\varphi_1$ in $T_1$. By the remarks
above, there exists a quantifier-free ${\cal L}_0$-formula $\varphi_0$ equivalent to $\varphi_1$ in $T_1$.
\end{proof}

\begin{cor}\label{PrimCor2}
There exists a quantifier-free ${\cal L}_0$-formula $\psi_0(C)$ such that for all perfect fields $K$ and
$c\in K^M$,
$$K\models\psi_0(c)\quad\Longleftrightarrow\quad 
\text{$\bigl(f_1(c,X),\dots,f_n(c,X)\bigr)K[X]$ is prime.}$$
\end{cor}
\begin{proof}
By the previous corollary and  \eqref{RadicalEqu}, for a perfect field $K$ and $c\in K^M$.
\end{proof}

Let now ${\cal L}_2$ be the language of rings $\cal L$, augmented by function symbols $^{-1}$ (unary)
and $\lambda_{mi}$ ($m$-ary), for all $1\leq i\leq m$. We extend the theory of rings
to an ${\cal L}_2$-theory $T_2$ by adding the defining axiom \eqref{InvAxiom} and for each $m\geq 1$ an axiom saying that for any
model of $T_2$ with underlying field $K$ and $a=(a_1,\dots,a_m)\in K^m$, the vector $\lambda_m(a)=
\bigl(\lambda_{m1}(a),\dots,\lambda_{mm}(a)\bigr)\in K^m$ is a non-trivial solution of the 
equation
$$a_1Y_1^p + \cdots + a_mY_m^p=0,$$
if there is such a solution and $\operatorname{char} K = p>0$. Every field may be 
expanded to a model of $T_2$.
Note that $T_2$ is a {\em universal}\/
theory, and if $K\subseteq L$ are the underlying fields of an extension of models of $T_2$,
then $L|K$ is a separable field extension. Along the lines of the proof of Corollary~\ref{PrimaryCor},
using part (1) of Corollary~\ref{ExtensionIdeal} instead of (2), one shows:

\begin{cor}\label{PrimCor3}
There exists a quantifier-free ${\cal L}_2$-formula $\varphi_2(C)$ such that for every field $K$ and all
$c\in K^M$,
$$K\models\varphi_2(c)\quad\Longleftrightarrow\quad 
\text{$\bigl(f_1(c,X),\dots,f_n(c,X)\bigr)K[X]$ is radical.}$$
\end{cor}

Hence in particular, the quantifier-free formula $\psi_2=\varphi_1\wedge\varphi_2$ in the language
${\cal L}_1\cup{\cal L}_2$ defines the set \eqref{PrimalitySet} for all fields $R=K$.

\subsection*{Prime ideals in polynomial rings over some arithmetical rings}

Based on the previous results, it is now straightforward
to produce numerous corollaries about the definability of 
primality for ideals in
polynomial rings $R[X]$, where $R$ is a DVR, a PID, etc. 
In order to keep the notational effort minimal, we restrict
ourselves to treating the following two situations:
\begin{enumerate}
\item $R$ is a DVR with perfect residue and fraction field.
\item $R=\Z$.
\end{enumerate}
The following lemma, whose proof is left to the reader, is fundamental.
Let $R$ be a domain, $f_1,\dots,f_n\in R[X]$, and $I=(f_1,\dots,f_n)R[X]$.
By Hermann's Theorem (Corollary~\ref{Artinian-Cor-2} for $R$ local and $l=1$) 
and Cramer's Rule, there exists a non-zero 
$\delta=\delta(f_1,\dots,f_n)\in R$ such that $IF[X]\cap R[X]=I:\delta R[X]$.

\begin{lemma}\label{PrimeCrit}
Suppose that $R$ is a B\'ezout domain. Then $I$ is a prime ideal
if and only if one of the following holds:
\begin{enumerate}
\item $IF[X]$ is a prime ideal of $F[X]$ and $I:\delta R[X]=I$, or
\item there exists a prime factor $t$ of $\delta$ such that
$t\in I$, and the image of $I$ in $(R/t R)[X]$ is a prime ideal.
\end{enumerate}
\end{lemma}

Let ${\cal L}_{\div}^*$ be the language
$${\cal L}_{\div}^* = {\cal L}_{\div} \cup \{Z_m : m\geq 1\} \cup  \{\bar{Z}_m : m\geq 1\},$$
where $Z_m$ and $\bar{Z}_m$ are $m$-ary predicate symbols, for $m\geq 1$. 
We construe a valuation ring $R$ (with fraction field $F$ and residue
field $\bar{R}$) 
as an ${\cal L}_{\div}^*$-structure as follows: we interpret the
symbols of ${\cal L}_{\div}$ as usual, and for $a_1,\dots,a_{m}\in R$, $m\geq 1$, 
we put
$$
R\models Z_m(a_1,\dots,a_{m}) \ :\Longleftrightarrow\
\text{$T^m+a_{1}T^{m-1}+\cdots+a_m \in R[T]$ %is separable and 
has a zero in $F$}
$$
and
$$
R\models\bar{Z}_m(a_1,\dots,a_{m}) \ :\Longleftrightarrow\ 
\text{$T^m+\bar{a_{1}}T^{m-1}+\cdots+\bar{a_m} \in \bar{R}[T]$ %is separable and 
has a zero in $\bar{R}$.}
$$
From Corollaries~\ref{DefBezout-Cor-3}, 
\ref{SolvDVR}, \ref{PrimCor2} and Lemma~\ref{PrimeCrit},
we get:

\begin{cor} 
There exists
a quantifier-free ${\cal L}_{\div}^*$-formula $\pi(C,T)$ 
such that for all DVRs $R$ with maximal ideal $\frak m$,
perfect residue field $\bar{R}=R/\frak m$ and perfect fraction field $F$, all 
generators $t$ of $\m$ and all $c\in R^M$,
$$R\models\pi(c,t)
\qquad\Longleftrightarrow\qquad \text{$\bigl(f_1(c,X),\dots,f_n(c,X)\bigr)R[X]$
is a prime ideal.}$$
\end{cor} 

The corollary above applies in particular to the ring $R=\Z_p$ 
of $p$-adic integers (with finite residue field $\F_p$ and
fraction field $\Q_p$ of characteristic zero). Let ${\cal L}_{\pow}$ be the
language obtained by augmenting the language ${\cal L}$ of rings by a unary predicate symbols $P_n$,
for each $n>0$. We construe $\Z_p$ as an ${\cal L}_{\pow}$-structure by interpreting the ring symbols
as usual and $P_n$ by the set
$$\bigl\{a\in\Z_p : \exists b\in\Z_p: b^n=a\bigr\}.$$
By Macintyre's Theorem \cite{Macintyre}, the complete ${\cal L}_{\pow}$-theory of $\Z_p$ admits
quan\-ti\-fier-elimination. Clearly, the relations on $\Z_p^m$ given by
$Z_m$ and $\bar{Z}_m$ are definable in the
${\cal L}_{\div}$-structure $\Z_p$. Hence:

\begin{cor}\label{DefPrime}
For each prime $p$, there exists a quantifier-free ${\cal L}_{\pow}$-formula $\pi_p(C)$ such that for all $c\in\Z_p^M$,
$$\Z_p\models\pi_p(c)\qquad\Longleftrightarrow\qquad \text{$\bigl(f_1(c,X),\dots,f_n(c,X)\bigr)\Z_p[X]$
is a prime ideal.}$$  
\end{cor}

\begin{example}
Let $R$ be a valuation ring and
suppose that $t$ is an element of $R$ of smallest
positive valuation. For $c\in R$ consider the ideal
$I_{c}=\bigl(t(1-tX),t cX\bigr)$ of $R[X]$, where $X$ is a single
indeterminate. Then we have 
$$\text{$I_{c}$ is a prime ideal} \quad\Longleftrightarrow\quad
\text{$t\in\sqrt{cR}$.}$$
To see this, use Lemma~\ref{PrimeCrit} and note that 
by the remark following Corollary~\ref{CorRadical}
we have $t\in I_{c}$ if and only if $t\in\sqrt{cR}$. 
Letting $R$ range over all $p$-adically closed valuation rings
(= models of the complete ${\cal L}_{\pow}$-theory of $\Z_p$) 
and taking $t=p$, this implies:
{\it there exists no ${\cal L}_{\pow}$-formula $\varphi(C)$ with the property
that for all $p$-adically closed valuation rings and all $c\in R$, we have
$R\models\varphi(c)$ if and only if $I_c$ is a prime ideal of $R[X]$.}
\end{example}

For homogeneous ideals, however, we do have a more uniform version of 
Corollary~\ref{DefPrime}:

\begin{prop}\label{Homogeneous}
Suppose that $f_1(C,X),\dots,f_n(C,X)\in\Z[C,X]$ are homogeneous in
$X=(X_1,\dots,X_N)$.
Then there exists a quantifier-free ${\cal L}_{\pow}$-formula $\pi_p'(C)$
such that for all $p$-adically closed 
valuation rings $R$ and all $c\in R^M$:
$$R\models\pi_p'(c)
\qquad\Longleftrightarrow\qquad \text{$\bigl(f_1(c,X),\dots,f_n(c,X)\bigr)R[X]$
is a prime ideal.}$$
\end{prop}
\begin{proof}
This follows from the remark following Corollary~\ref{SolvDVR},
Corollaries~\ref{DefBezout-Cor-3}, \ref{PrimCor2}, and Lemma~\ref{PrimeCrit}.
\end{proof}

We now extend ${\cal L}_{\rad}$ to a language ${\cal L}_{\rad}^*$ by adjoining, for each $m\geq 1$, an
$(m+1)$-ary predicate symbol $Z_m$. We expand the ${\cal L}_{\rad}$-structure $\Z$ to an 
${\cal L}_{\rad}^*$-structure by interpreting the $Z_m$ as follows: for $p,a_1,\dots,a_{m}\in\Z$,
\begin{multline*}
\Z\models Z_m(p,a_1,\dots,a_{m}) \qquad\Longleftrightarrow\qquad\\
\text{$p$ is a prime and 
$T^m+\bar{a_{1}}T^{m-1}+\cdots+\bar{a_m} \in \F_p[T]$ has a zero in $\F_p$, or} \\
\text{$p=0$ and 
$T^m+a_{1}T^{m-1}+\cdots+a_m \in \Z[T]$ has a zero in $\Q$.}
\end{multline*}
Let us call an ${\cal L}_{\rad}^*$-formula $\varphi(C)$ {\bf special} if it is of the form
$$\varphi(C)=\exists U\bigl(\text{``$U$ is a prime factor of $\delta(C)$''} \wedge \psi(C,U)\bigr),$$
where $U$ is a single new variable, $\delta(C)\in\Z[C]$,
and $\psi(C,U)$ a quantifier-free ${\cal L}_{\rad}^*$-formula. Using
the remark following Corollary~\ref{DefBezout-Cor-3}, 
Theorem~\ref{SolvBezout}, Lemma~\ref{PrimeCrit} and
Corollary~\ref{PrimCor2}, we get:
 
\begin{cor}
There exists a finite disjunction $\pi(C)$ of special ${\cal L}_{\rad}^*$-for\-mu\-las such that for
all $c\in\Z^M$,
$$\Z\models\pi(c) 
\qquad\Longleftrightarrow\qquad \text{$\bigl(f_1(c,X),\dots,f_n(c,X)\bigr)\Z[X]$
is a prime ideal.}$$ 
\end{cor}

We leave it to the reader to formulate a result analogous to
Proposition~\ref{Homogeneous} for homogeneous ideals in polynomial rings over 
B\'ezout domains.

\section{Appendix}\label{Appendix}

We would like to point out another application of the useful
Proposition~\ref{NilProp}, to a
characterization of Jacobson domains among Noetherian domains.
A {\bf Jacobson ring} is a ring each of whose prime ideals is an
intersection of maximal ideals. (See \cite{Bourbaki}, IV.3.4 for the
basic properties stated below.) 
A Jacobson domain is a Jacobson ring which happens to be a
domain. Examples for Jacobson domains include $\Z$ (or more generally: any 
PID with infinitely many pairwise non-associated primes), and 
every polynomial ring $R[X]$ over a Jacobson domain $R$.
For a domain $R$, we denote the integral closure of $R$ in an 
algebraic closure of $\Frac(R)$ by $\tilde{R}$.

\begin{prop}\label{PS-Prop}
Let $R$ be a Jacobson domain and $f_0,\dots,f_n\in R[X]$. Then
$$f_0\in\sqrt{(f_1,\dots,f_n)R[X]} \quad\Longleftrightarrow\quad
\text{$f_0(a)\in\sqrt{\bigl(f_1(a),\dots,f_n(a)\bigr)\tilde{R}}$ for all $a\in\tilde{R}^N$.}$$
\end{prop}
\begin{proof}
The implication $\Rightarrow$ is clear. To prove $\Leftarrow$, assume $f_0\notin\sqrt{I}$,
where $I=(f_1,\dots,f_n)R[X]$. 
We have to find $a\in\tilde{R}^N$ such that
$$f_0(a)\notin\sqrt{\bigl(f_1(a),\dots,f_n(a)\bigr)\tilde{R}}.$$ 
We write $f_i(X)$ as $f_i(c,X)$, with $f_i(C,X)\in\Z[C,X]$, $c\in R^M$, for
$i=0,\dots,n$.
By Proposition~\ref{NilProp}, there exists
$i\in\{1,\dots,m\}$ with $$q_i(c)\notin\sqrt{\bigl(p_{i1}(c),\dots,p_{ik}(c)\bigr)R}.$$
Take a maximal
ideal $\m$ of $R$ that contains $\sqrt{\bigl(p_{i1}(c),\dots,p_{ik}(c)\bigr)R}$ but not $q_i(c)$, and a
maximal ideal $\tilde{\m}$ of $\tilde{R}$ lying above $\m$. Then by \eqref{Nil1} and
\eqref{Nil2} for the algebraically closed field $K=\tilde{R}/\tilde{\m}$, there exists $a\in\tilde{R}^N$
with $$f_0(c,a)\notin\tilde{\m}, f_1(c,a),\dots,f_n(c,a)\in\tilde{\m}.$$ 
Hence
$f_0(c,a)\notin\sqrt{\bigl(f_1(c,a),\dots,f_n(c,a)\bigr)\tilde{R}}$. 
\end{proof}

\begin{remark}
The case $R=\Z$ of the last proposition is Theorem~5.3 in 
\cite{Prestel-Schmid}. (The proof in \cite{Prestel-Schmid} is much longer.)
\end{remark}

\begin{cor}\label{Jacobson-Char}
Let $R$ be a Noetherian domain. The following are equivalent:
\begin{enumerate}
\item $R$ is a Jacobson domain.
%\item For all $N\in\N$ and  
%$f_0,f_1,\dots,f_n\in R[X_1,\dots,X_N]$, if
%$$f_0(a)\in\Jac\bigl(f_1(a),\dots,f_n(a)\bigr)\quad\text{for 
%all $a\in\tilde{R}^N$,}$$ then $$f_0\in \sqrt{(f_1,\dots,f_n)R[X]}.$$
\item For all $N\in\N$ and  
$f_0,f_1,\dots,f_n\in R[X_1,\dots,X_N]$, if
$$f_0(a)\in\sqrt{\bigl(f_1(a),\dots,f_n(a)\bigr)\tilde{R}}\quad\text{for 
all $a\in\tilde{R}^N$,}$$
then $$f_0\in \sqrt{(f_1,\dots,f_n)R[X]}.$$
%\item For all $N\in\N$ and
%$f_0,f_1,\dots,f_n\in R[X_1,\dots,X_N]$, $$f_0(a)\in
%\bigl(f_1(a),\dots,f_n(a)\bigr)\text{\ for all $a\in\tilde{R}^N$}
%\ \Rightarrow\ f_0\in \sqrt{(f_1,\dots,f_n)}.$$
\item For polynomials
$f_1,\dots,f_n\in R[X]$ in the single indeterminate $X$, if
$$1\in\bigl(f_1(a),\dots,f_n(a)\bigr)\tilde{R}\quad
\text{for all $a\in\tilde{R}$}$$
then
$$1\in(f_1,\dots,f_n)R[X].$$
\end{enumerate}
\end{cor}
\begin{proof}
The implications (1) $\Rightarrow$ (2) $\Rightarrow$ (3)
don't need the assumption that
$R$ be Noetherian:
(1)~$\Rightarrow$~(2) follows from the proposition,
and (2)~$\Rightarrow$~(3) is trivial.
Assume now that $R$ is a Noetherian domain, and (3) holds. 
In order to show that $R$ is Jacobson 
it then suffices to show the following: if $r_1,\dots,r_n,r\in R$
have the property that $r\in\frak m$ for every maximal ideal
$\frak m$ of $R$ which contains $r_1,\dots,r_n$, then 
$r\in\sqrt{(r_1,\dots,r_n)R}$.
For this we may assume $r\neq 0$, and we 
consider the polynomials $r_1,\dots,r_n,1-rX$ in the single
indeterminate $X$. We claim that for every $a\in\tilde{R}$, we have
$1\in (r_1,\dots,r_n,1-ra)\tilde{R}$.
Suppose otherwise, and let $a\in\tilde{R}$ with
$1\notin (r_1,\dots,r_n,1-ra)\tilde{R}$. 
Let $\frak n$ be a maximal ideal of $\tilde{R}$
containing $r_1,\dots,r_n,1-ra$. Then
${\frak m}:={\frak n}\cap R$ is a maximal ideal of $R$ (by the going-up
property for integral ring extensions).
Since $r_1,\dots,r_n\in{\frak m}$ we get $r\in\frak m$ and
hence $1\in\frak n$, a contradiction. By (3), this implies that there
exists a relation
$$1=r_1g_1+\cdots+r_ng_n+(1-rX)g,$$
where $g_1,\dots,g_n,g\in R[X]$. Substituting $1/r$ for $X$ and 
multiplying both sides by $r^d$, where $d$ is the maximum of the degrees of
the $g_i$, we get $r^d \in (r_1,\dots,r_n)R$. (Rabinowitsch trick.) Hence
$r\in\sqrt{(r_1,\dots,r_n)R}$ as desired.
\end{proof}

The following easily proved
lemma (applied to $R=S=\tilde{\Z}$) gives an example which shows that
in Proposition~\ref{PS-Prop} above the $\sqrt{\cdots}$ cannot be omitted:

\begin{lemma}\label{Counterexample}
Let $R$ be a Pr\"ufer domain and $$f_0(X,Y)=X^2,\ f_1(X,Y)=X^2+Y^2,\ 
f_2(X,Y)=XY,$$ where $X$ and $Y$ are single indeterminates. 
Then $f_0(x,y)\in\bigl(f_1(x,y),f_2(x,y)\bigr)R$ for all
$(x,y)\in R^2$, but $f_0\notin (f_1,f_2)S[X,Y]$ for every domain $S$ extending $R$.
\end{lemma}

For the polynomials $f_1=-2,f_2=X^2+X+1\in\Z[X]$ (where $X$ is a
single indeterminate) we have $1\in\bigl(f_1(a),f_2(a)\bigr)\Z$ for every
$a\in\Z$, but $1\notin(f_1,f_2)\Z[X]$. This (well-known) example shows that in
Proposition~\ref{PS-Prop} we cannot replace $\tilde{R}$ by $R$.
On the other hand, we do have
$1=gf_1+f_2$ where $g=\frac{X(X+1)}{2}$ is
integer-valued, that is, $g(\Z)\subseteq\Z$. Indeed, Skolem \cite{Skolem}
showed in general:

\begin{prop}
If $f_0,f_1,\dots,f_n\in\Z[X]=
\Z[X_1,\dots,X_N]$
satisfy $$f_0(a)\in\sqrt{\bigl(f_1(a),\dots,f_n(a)\bigr)\Z}\quad
\text{for all $a\in\Z^N$,}$$
then $$f_0\in\sqrt{(f_1,\dots,f_n)\operatorname{Int}(\Z^N)},$$ where 
$$\operatorname{Int}(\Z^N)=\left\{f(X)\in\Q[X]: \text{$f(a)\in\Z$ 
for all $a\in\Z^N$}\right\}$$ denotes the subring
of $\Q[X]$ of integer-valued polynomials.
\end{prop}
One says that
the domain $\Z$ has the {\bf Skolem property.}
See \cite{Frisch} for a characterization of the Noetherian domains with
the Skolem property similar in spirit to Corollary~\ref{Jacobson-Char}.
In \cite{Brizolis} it is shown that for $N=1$, 
the $\sqrt{\cdots}$ in Skolem's theorem may be omitted.
Lemma~\ref{Counterexample} above (for $R=\Z$, $S=\Q$)
shows that for $N>1$ we cannot omit the
$\sqrt{\cdots}$. 
This gives a negative answer to a question posed in
\cite{Brizolis}.

\bibliographystyle{amsplain}
\bibliography{ideal3.bib}

\providecommand{\bysame}{\leavevmode\hbox to3em{\hrulefill}\thinspace}
\providecommand{\MR}{\relax\ifhmode\unskip\space\fi MR }
\providecommand{\MRhref}[2]{%
  \href{http://www.ams.org/mathscinet-getitem?mr=#1}{#2}
}
\providecommand{\href}[2]{#2}
\begin{thebibliography}{10}

\bibitem{maschenb-thesis}
M.~Aschenbrenner, \emph{Ideal {M}embership in {P}olynomial {R}ings over the
  {I}ntegers}, Ph.D. thesis, University of Illinois at Urbana-Champaign, 2001.

\bibitem{maschenb-ideal2}
\bysame, \emph{Ideal membership in polynomial rings over the integers},
  submitted, 2003.

\bibitem{Ax}
J.~Ax, \emph{Solving diophantine problems modulo every prime}, Ann. of Math.
  (2) \textbf{85} (1967), 161--183.

\bibitem{Bass}
H.~Bass, \emph{Torsion free and projective modules}, Trans. Amer. Math. Soc.
  \textbf{102} (1962), 319--327.

\bibitem{Bourbaki}
N.~Bourbaki, \emph{{\'E}l\'ements de {M}ath\'ematique. {A}lg\`ebre
  {C}ommutative}, Hermann, Paris, 1964.

\bibitem{BA}
\bysame, \emph{{\'El\'ements de Math\'ematique. Alg\`ebre. Chapitres 4 \`a 7}},
  Lecture Notes in Mathematics, vol. 864, Masson, Paris, 1981.

\bibitem{Brizolis}
D.~Brizolis, \emph{A theorem on ideals in {P}r\"ufer rings of integral-valued
  polynomials}, Comm. Algebra \textbf{7} (1979), no.~10, 1065--1077.

\bibitem{Chang-Keisler}
C.~C. Chang and H.~J. Keisler, \emph{Model {T}heory}, Studies in Logic and the
  Foundations of Mathematics, vol.~73, North-Holland Publishing Co., Amsterdam,
  1973.

\bibitem{Chase}
S.~U. Chase, \emph{Direct products of modules}, Trans. Amer. Math. Soc.
  \textbf{97} (1960), 457--473.

\bibitem{Cherlin}
G.~Cherlin, \emph{Algebraically closed commutative rings}, J. Symbolic Logic
  \textbf{38} (1973), 493--499.

\bibitem{Cohen}
I.~S. Cohen, \emph{Commutative rings with restricted minimum condition}, Duke
  Math J. \textbf{17} (1950), 27--42.

\bibitem{vdDries-Thesis}
L.~van~den Dries, \emph{Model {T}heory of {F}ields}, Ph.D. thesis, R.U.
  Utrecht, 1978.

\bibitem{vdDries-Holly}
L.~van~den Dries and J.~Holly, \emph{Quantifier elimination for modules with
  scalar variables}, Ann. Pure Appl. Logic \textbf{57} (1992), no.~2, 161--179.

\bibitem{vdDries-Schmidt}
L.~van~den Dries and K.~Schmidt, \emph{Bounds in the theory of polynomial rings
  over fields. {A} nonstandard approach}, Invent. Math. \textbf{76} (1984),
  no.~1, 77--91.

\bibitem{Eisenbud-Evans}
D.~Eisenbud and E.~G. Evans, \emph{Generating modules efficiently: theorems
  from algebraic {K}-theory}, J. Algebra \textbf{27} (1973), 278--305.

\bibitem{Forster}
O.~Forster, \emph{{\"U}ber die {A}nzahl der {E}rzeugenden eines {I}deals in
  einem noetherschen {R}ing}, Math. Z. \textbf{84} (1964), 80--87.

\bibitem{Frisch}
S.~Frisch, \emph{Nullstellensatz and {S}kolem properties for integer-valued
  polynomials}, J. Reine Angew. Math. \textbf{536} (2001), 31--42.

\bibitem{Glaz}
S.~Glaz, \emph{Commutative {C}oherent {R}ings}, Lecture Notes in Math., vol.
  1371, Springer-Verlag, Berlin-Heidelberg-New York, 1989.

\bibitem{Glaz-History}
\bysame, \emph{Commutative coherent rings: historical perspective and current
  developments}, Nieuw Arch. Wisk. \rom{(4)} \textbf{10} (1992), no.~1--2,
  37--56.

\bibitem{Greenberg-Vasconcelos}
B.~Greenberg and W.~Vasconcelos, \emph{Coherence of polynomial rings}, Proc.
  Amer. Math. Soc. \textbf{54} (1976), 59--64.

\bibitem{Heitmann-1}
R.~Heitmann, \emph{Generating ideals in {P}r\"ufer domains}, Pacific J. Math.
  \textbf{62} (1976), no.~1, 117--126.

\bibitem{Heitmann-2}
\bysame, \emph{Generating non-noetherian modules efficiently}, Michigan Math.
  J. \textbf{31} (1984), 167--180.

\bibitem{Hermann}
G.~Hermann, \emph{{Die Frage der endlich vielen Schritte in der Theorie der
  Polynom\-ideale}}, Math. Ann. \textbf{95} (1926), 736--788.

\bibitem{Hodges}
W.~Hodges, \emph{Model {T}heory}, Encyclopedia of Mathematics and its
  Applications, vol.~42, Cambridge University Press, Cambridge, 1993.

\bibitem{Kollar}
J.~Koll\'ar, \emph{Sharp effective {N}ullstellensatz}, J. Amer. Math. Soc.
  \textbf{1} (1988), 963--975.

\bibitem{Macintyre}
A.~Macintyre, \emph{On definable subsets of $p$-adic fields}, J. Symbolic Logic
  \textbf{41} (1976), 605--610.

\bibitem{Northcott}
D.~G. Northcott, \emph{A generalization of a theorem on the content of
  polynomials}, Proc. Cambridge Philos. Soc. \textbf{55} (1959), 282--288.

\bibitem{OLeary-Vaaler}
R.~O'Leary and J.~Vaaler, \emph{Small solutions to inhomogeneous linear
  equations over number fields}, Trans. Amer. Math. Soc. \textbf{336} (1993),
  no.~2, 915--931.

\bibitem{Prestel-Schmid}
A.~Prestel and J.~Schmid, \emph{Existentially closed domains with radical
  relations}, J. Reine Angew. Math. \textbf{407} (1990), 178--201.

\bibitem{Sabbagh}
G.~Sabbagh, \emph{Coherence of polynomial rings and bounds in polynomial
  ideals}, J. Algebra \textbf{31} (1974), 499--507.

\bibitem{Schoutens}
H.~Schoutens, \emph{Bounds in polynomial rings over {A}rtinian local rings},
  manuscript, 2002.

\bibitem{Schuelting}
H.-W. Sch{\"u}lting, \emph{\"{U}ber die {E}rzeugendenanzahl invertierbarer
  {I}deale in {P}r\"uferringen}, Comm. Algebra \textbf{7} (1979), no.~13,
  1331--1349.

\bibitem{Seidenberg1}
A.~Seidenberg, \emph{Constructions in algebra}, Trans. Amer. Math. Soc.
  \textbf{197} (1974), 273--313.

\bibitem{Skolem}
T.~Skolem, \emph{Ein {S}atz \"uber ganzwertige {P}olynome}, Norske Vid. Selsk.
  Forh. \textbf{9} (1936), 111--113.

\bibitem{Soublin-2}
J.-P. Soublin, \emph{Un anneau coh\'erent dont l'anneau des polynomes n'est pas
  coh\'erent}, C. R. Acad. Sci. Paris S\'er. A \textbf{267} (1968), 241--243.

\bibitem{Soublin-1}
\bysame, \emph{Anneaux et modules coh\'erents}, J. Algebra \textbf{15} (1970),
  455--472.

\bibitem{Swan-1}
R.~Swan, \emph{The number of generators of a module}, Math. Z. \textbf{102}
  (1967), 318--322.

\bibitem{Swan-2}
\bysame, \emph{{$n$}-generator ideals in {P}r\"ufer domains}, Pacific J. Math.
  \textbf{111} (1984), no.~2, 433--446.

\end{thebibliography}

\end{document}